\documentclass[a4paper]{amsart}
\usepackage{amssymb} 
\usepackage[T1]{fontenc}
\usepackage[latin1]{inputenc}
\usepackage{amsfonts}
\usepackage{amsxtra}
\usepackage{ae}
\usepackage[all]{xy}
\usepackage{enumerate}

\include{diagram}

\newcommand*{\Cov}{\mathfrak{Cov}}

\newcommand*{\D}{\mathcal{D}}

\newcommand*{\cotimes}{\hat{\otimes}}
\newcommand*{\SHom}{\mathfrak{Hom}}
\newcommand*{\SExt}{\mathfrak{Ext}}
\newcommand*{\LSMod}{\text{-}\mathfrak{Mod}}

\newcommand*{\Or}{\mathsf{Or}}

\DeclareMathOperator{\Sh}{Sh}
\DeclareMathOperator{\ind}{ind}
\DeclareMathOperator{\pro}{pro}

\DeclareMathOperator{\Tr}{tr}

\DeclareMathOperator{\Hom}{Hom}
\DeclareMathOperator{\map}{Map}
\DeclareMathOperator{\im}{im}
\DeclareMathOperator{\id}{id}

\DeclareMathOperator{\hyp}{hyp}
\DeclareMathOperator{\elli}{ell}

\newenvironment{bnum}
{\begin{list}{}
    {\setlength{\labelwidth}{15pt}
     \setlength{\leftmargin}{\labelwidth}
    }
}
{\end{list}}

\numberwithin{equation}{section}
\theoremstyle{change}
\newtheorem{theorem}{Theorem}[section]
\newtheorem{prop}[theorem]{Proposition}
\newtheorem{lemma}[theorem]{Lemma}
\newtheorem{cor}[theorem]{Corollary}
\newtheorem{definition}[theorem]{Definition}

\begin{document}

\title[Equivariant cohomology for totally disconnected groups]{A new description of equivariant cohomology for 
totally disconnected groups}
\author{Christian Voigt}
\address{Mathematisches Institut\\
         Westfälische Wilhelms-Universität Münster\\
         Einsteinstr.\ 62\\
         48149 Münster\\
         Germany
}
\email{cvoigt@math.uni-muenster.de}


\thanks{This research was supported by the EU-Network \emph{Quantum
    Spaces and Noncommutative Geometry} (Contract HPRN-CT-2002-00280)
  and the \emph{Deutsche Forschungsgemeinschaft} (SFB 478).}

\maketitle

\begin{abstract}
We consider smooth actions of totally disconnected groups on simplicial complexes and compare 
different equivariant cohomology groups associated to such actions. Our main result is that the 
bivariant equivariant cohomology theory introduced by Baum and Schneider can be described using 
equivariant periodic cyclic homology. 
This provides a new approach to the construction of Baum and Schneider as well 
as a computation of equivariant periodic cyclic homology for a natural class of examples.
In addition we discuss the relation between cosheaf homology and equivariant 
Bredon homology. 
Since the theory of Baum and Schneider generalizes cosheaf homology we finally see 
that all these approaches to equivariant cohomology for totally disconnected 
groups are closely related. 
\end{abstract}

\section{Introduction}

In this paper we compare different cohomology theories defined for smooth 
proper actions of totally disconnected groups on appropriate topological spaces. 
More precisely, we consider bivariant equivariant cohomology in the sense of  
Baum and Schneider \cite{BS}, equivariant periodic cyclic homology \cite{Voigtperiodic}, cosheaf 
homology \cite{BCH} and equivariant Bredon homology \cite{Bredon1}, \cite{Lueck}. Our main 
result explains the relation between the first two theories. Moreover, by the work of 
Baum and Schneider, cosheaf cohomology can be viewed as a special case of bivariant 
equivariant cohomology. We complete the 
picture by showing that equivariant Bredon homology is naturally isomorphic to cosheaf homology 
for proper actions on simplicial complexes. \\
The approach of Baum and Schneider is based on sheaf theory and unifies several 
constructions which appeared previously in the literature. As already mentioned, it contains 
as a special case the cosheaf homology groups considered by Baum, Connes and Higson 
in connection with the Baum-Connes conjecture \cite{BCH}. Moreover it covers a 
construction of Baum and Connes for discrete groups \cite{BC2}. 
Equivariant Bredon homology is an equivariant generalization of cellular homology. 
It has been used by L\"uck to describe the rationalized left 
hand side of the assembly maps in $ K $- and $ L $-theory for discrete groups \cite{Lueck}. 
Finally, equivariant cyclic homology can be viewed as a 
noncommutative generalization of the equivariant de Rham cohomology of manifolds. \\
The latter theory is different in nature to the previous ones since it is defined not only for spaces 
but also for possibly noncommutative algebras equipped with a group action. Its construction 
is based on an extension of the Cuntz-Quillen approach to cyclic homology. However, computations tend to be  
more difficult than in ordinary cyclic homology due to the fact that the basic ingredient in the theory is 
not a complex in the usual sense of homological algebra. In fact, our original motivation 
was to calculate equivariant periodic cyclic homology for some basic examples. \\
Let us now describe in more detail the contents of this paper. In section \ref{sectotdis} we review some 
facts about totally disconnected groups and smooth representations. We recall the concept 
of a covariant module which plays a central role in equivariant cyclic homology. In particular 
we discuss the decomposition of a covariant module into an 
elliptic part and a hyperbolic part in the case of totally disconnected groups. 
In section \ref{secper} we review the general construction of equivariant cyclic homology. Bivariant 
equivariant periodic cyclic homology $ HP^G_*(A,B) $ is defined for pairs of $ G $-algebras $ A $ and $ B $.  
We are interested in the case of certain algebras of smooth functions 
on simplicial complexes acted upon by a totally disconnected group $ G $. The appropriate notion of 
a smooth function on a simplicial complex $ X $ is introduced in section \ref{secsmooth} where we also study some 
properties of the resulting algebra $ C^\infty_c(X) $ of smooth functions with compact support. 
Let us point out that simplicial complexes are a convenient choice of spaces for two reasons. On the one hand 
they are special enough to have a nice de 
Rham-theoretic description of their cohomology. On the other hand they are general enough to cover 
natural examples, most notably in connection with the Baum-Connes conjecture. 
In section \ref{secBSdef} we review the definition of bivariant equivariant cohomology in the sense 
of Baum and Schneider. We study in particular the case of group actions on simplicial complexes and  
introduce the notion of a $ G $-simplicial complex. 
Section \ref{secHKR} contains an equivariant version of the Hochschild-Kostant-Rosenberg theorem. 
This theorem is an important ingredient in our main result which 
is formulated and proved in section \ref{secBSC}. We show that there exists a natural isomorphism 
$$
HP^G_*(C^\infty_c(X), C^\infty_c(Y)) \cong \bigoplus_{j \in \mathbb{Z}} H^{* + 2j}_G(X,Y) 
$$
if $ X $ and $ Y $ are finite dimensional locally finite $ G $-simplicial complexes and $ X $ is proper. 
Here $ H^*_G $ denotes the theory of Baum and Schneider.  
A small variant of this result in the case of discrete groups acting 
on manifolds yields a description of the theory of Baum and Connes \cite{BC2} in terms of 
equivariant cyclic homology. 
Finally, in section \ref{secbred} we review the definitions of equivariant Bredon homology and 
cosheaf homology and show that these theories are naturally isomorphic on the category of $ G $-simplicial complexes. \\
This paper is based on the last chapter of my thesis ~\cite{Voigtthesis}. I would like to thank J. Cuntz for his constant 
support and P. Schneider for some helpful comments. 

\section{Totally disconnected groups}\label{sectotdis}

Let $ G $ be locally compact and totally disconnected group. We call an element $ t \in G $ elliptic if it is 
contained in a compact subgroup. The set of all elliptic elements of $ G $ is denoted by $ G_{\elli} $. It is known that 
$ G_{\elli} $ is a closed subset of $ G $ \cite{Willis}. In contrast 
we shall say that an element $ t \in G $ is hyperbolic if it is not elliptic. Let $ G_{\hyp} $ be the set 
of all hyperbolic elements of $ G $. Hence, according to these definitions, we obtain a disjoint union decomposition 
$$ 
G = G_{\elli} \cup G_{\hyp} 
$$ 
of the space $ G $. Throughout the paper we shall assume that $ G_{\elli} $ is the union of all compact open subgroups of $ G $. 
In this case $ G_{\hyp} $ is again open and closed. We will refer to 
groups satisfying this condition simply as totally disconnected groups. If $ G $ acts properly on an affine 
Bruhat-Tits building then it is a totally disconnected group in this sense ~\cite{HN}. Other examples are of course all discrete groups. \\
In connection with equivariant cyclic homology we work with smooth representations of totally disconnected groups 
on bornological vector spaces. For general information on bornological vector spaces and smooth representations of locally compact 
groups we refer to \cite{Meyerthesis}, \cite{Meyersmoothrep}. \\
Let us recall some basic facts concerning smooth representations and fix our notation. 
All bornological vector spaces in this paper are assumed to be convex and complete.
A representation of a totally disconnected group $ G $ on a bornological vector space $ V $ is 
called smooth if the the stabilizers of small subsets in $ V $ are open subgroups of $ G $. To avoid confusion 
we point out that $ t \in G $ is in the stabilizer of a subset $ S \subset V $ iff $ t \cdot v = v $ for 
all $ v \in S $. For instance, the trivial representation of $ G $ on any bornological vector space is a smooth representation. 
We will frequently also speak of $ G $-modules instead of smooth representations. A bounded linear map $ f: M \rightarrow N $ 
between $ G $-modules is called equivariant if it commutes with the action of $ G $. \\
If $ G $ is a totally disconnected group we denote by $ \D(G) $ the space of locally constant 
functions on $ G $ with compact support. This space is equipped with the fine bornology and elements of $ \D(G) $ are also referred 
to as smooth functions on $ G $ with compact support. The left regular representation of $ G $ on $ \D(G) $ 
given by 
$$
(s \cdot f)(t) = f(s^{-1}t)
$$
is a basic example of a smooth representation. 
\begin{lemma}
Let $ G $ be a totally disconnected group and let $ t \in G $ be hyperbolic. Then there are no fixed points 
in the left regular representation $ \D(G) $ for the action of $ t $. 
\end{lemma}
\proof Assume $ f \in \D(G) $ is a nonzero element satisfying $ t \cdot f = f $. It suffices to consider the case 
that $ f $ is the characteristic function of some compact open subset $ K $ of $ G $. In this 
case $ t \cdot f = f $ just means $ t \cdot K = K $. Since $ K $ is compact there exists a natural number $ n $ 
such that $ t^n \in H $ for some compact open subgroup $ H $ of $ G $. We deduce that the closed subgroup generated by $ t $ is compact. 
This contradicts the assumption that $ t $ is hyperbolic. \qed \\
A bornological algebra is a bornological vector space $ A $ with an associative multiplication given as a 
bounded linear map $ A \cotimes A \rightarrow A $ where $ \cotimes $ denotes the completed bornological tensor product. 
Remark that we do not require the existence of a unit in a bornological algebra. A basic example of a bornological 
algebra is the Hecke algebra of a totally disconnected group $ G $. It is obtained by equipping the space $ \D(G) $ 
with the convolution product 
$$
(f * g)(t) = \int f(s) g(s^{-1} t) ds 
$$
where $ ds $ is a fixed left Haar measure on $ G $. This algebra is unital iff the group $ G $ is discrete. 
We will denote the Hecke algebra of $ G $ again by $ \D(G) $. \\
A module $ M $ over a bornological algebra $ A $ is called nondegenerate if the module action $ A \cotimes M \rightarrow M $ is a 
bornological quotient map. We remark that the category of smooth representations of $ G $ is isomorphic to the 
category of nondegenerate $ \D(G) $-modules. \\
A $ G $-algebra is a bornological algebra $ A $ which is at the same time a $ G $-module such that the 
multiplication $ A \cotimes A \rightarrow A $ is equivariant. Here the tensor product $ A \cotimes A $ is equipped 
with the diagonal action as usual. A particular example of a $ G $-algebra is the 
algebra $ \mathcal{K}_G $ which is defined as follows. 
As a bornological vector space we have $ \mathcal{K}_G = \D(G) \cotimes \D(G) = \D(G \times G) $. 
The multiplication in $ \mathcal{K}_G $ is given by  
\begin{equation*}
(k \cdot l)(s,t) = \int_G k(s,r) l(r,t) dr
\end{equation*}
and the $ G $-action is defined by
\begin{equation*}
(r\cdot k)(s,t) = k(r^{-1}s, r^{-1}t).
\end{equation*}
This algebra can be viewed as a dense subalgebra of the algebra of compact operators $ \mathbb{K}(L^2(G)) $ on the 
Hilbert space $ L^2(G) $. \\
Next we recall the definition of a covariant module. For more details we refer to \cite{Voigtperiodic}. 
Let $ \mathcal{O}_G $ be the space $ \D(G) $ equipped with pointwise multiplication and the action of $ G $ by 
conjugation. A covariant module $ M $ is a smooth representation of 
$ G $ which is at the same time a nondegenerate $ \mathcal{O}_G $-module. 
The $ G $-module structure and the $ \mathcal{O}_G $-module structure 
are required to be compatible in the sense that 
$$
s \cdot (f \cdot m) = (s \cdot f) \cdot (f \cdot m)
$$
for all $ s \in G, f \in \mathcal{O}_G $ and $ m \in M $. A bounded linear map $ f: M \rightarrow N $ between covariant modules 
is called covariant if it is $ \mathcal{O}_G $-linear and equivariant. The category of covariant modules is denoted 
by $ G \LSMod $. \\
Covariant modules are closely related to 
equivariant sheaves over the space $ G $-space $ G $ with the adjoint action. 
\begin{prop} \label{eqsh}
Let $ G $ be a totally disconnected group. The functor of taking global sections with compact support induces an 
equivalence between the category $ \Sh_G(G) $ of equivariant sheaves over the $ G $-space $ G $ with the adjoint action  
and the category of $ G $-covariant modules whose underlying bornology is fine. 
\end{prop}
In other words, covariant modules can be viewed as equivariant sheaves except that we include a bornology as 
extra information. The proof of of proposition \ref{eqsh} is sketched in \cite{BS}. \\
The category of covariant modules is isomorphic to the category of nondegenerate $ \Cov(G) $-modules where 
$ \Cov(G) = \mathcal{O}_G \rtimes G $ is the smooth crossed product of $ \mathcal{O}_G $ with respect to the adjoint action. 
In particular $ \Cov(G) $ itself is a covariant module in a natural way. \\
Every covariant module $ M $ is equipped with a natural automorphism $ T: M \rightarrow M $. For $ M = \Cov(G) $ 
this automorphism is defined by 
$$
T(f)(s,t) = f(s,st)
$$
where we view an element of $ \Cov(G) $ as a smooth function with compact support on $ G \times G $, the first variable corresponding 
to $ \mathcal{O}_G $. To define the operator $ T $ for an arbitrary covariant module $ M $ one uses the canonical 
isomorphism $ \Cov(G) \cotimes_{\Cov(G)} M \cong M $ and applies the map $ T: \Cov(G) \rightarrow \Cov(G) $ from above 
to the first tensor factor in $ \Cov(G) \cotimes_{\Cov(G)} M $. \\
In some situations we will have to look at the stalks of the sheaf underlying a covariant module. Let $ t \in G $ be a point 
and consider the ideal $ I_t $ in $ \mathcal{O}_G $ of functions vanishing at $ t $. Clearly $ I_t $ is a prime ideal and hence we 
may consider the localisation $ M_t $ of a covariant module $ M $ at $ I_t $. Since $ G $ is totally disconnected 
the space $ M_t $ can be identified with $ M/I_tM $ in a natural way. The localisation $ M_t $ is no longer 
a covariant module in general, in addition to the natural $ \mathcal{O}_G $-module structure we only have an 
action of the centralizer $ G_t $ of the element $ t $ on $ M_t $. 
A basic observation is that a sequence $ 0 \rightarrow K \rightarrow E \rightarrow Q \rightarrow 0 $ of (fine) covariant modules is exact 
iff the localized sequences $ 0 \rightarrow K_t \rightarrow E_t \rightarrow Q_t \rightarrow 0 $ are exact 
for all $ t \in G $. \\
Since the group $ G $ is a disjoint union of the open sets of elliptic and hyperbolic elements we can define   
two multipliers $ P_{\elli} $ and $ P_{\hyp} $ of the algebra $ \mathcal{O}_G $ as follows. The multiplier $ P_{\elli} $ is 
the characteristic function of the set $ G_{\elli} $ whereas $ P_{\hyp} = 1 - P_{\elli} $ is the characteristic 
function of $ G_{\hyp} $. For every covariant module $ M $ we obtain a natural direct sum decomposition of covariant 
modules $ M = M_{\elli} \oplus M_{\hyp} $ where 
$ M_{\elli} = P_{\elli} \cdot M $ and $ M_{\hyp} = P_{\hyp} \cdot M $. \\
We conclude this section with the construction of a canonical projection on the elliptic part $ M_{\elli} $ 
of an arbitrary covariant module $ M $. 
Consider the operator $ T: M_{\elli} \rightarrow M_{\elli} $ and let $ m \in M_{\elli} $. 
By the definition of $ G_{\elli} $ and since $ M_{\elli} $ is a smooth representation there exists a 
natural number $ n $ such that 
$ T^n(m) = m $. We can thus write $ M_{\elli} $ as direct limit of the spaces 
$ M(n) $ consisting of all elements $ m $ such that $ T^n(m) = m $. Here the direct limit is taken over 
the natural numbers where by definition $ k \leq l $ iff $ k $ divides $ l $. 
We define a covariant map $ E: M(n) \rightarrow M $ by 
$$
E(m) = \frac{1}{n} \sum_{j = 0}^{n - 1} T^j(m).
$$ 
This definition is compatible with the structure maps in the inductive limit and yields 
a covariant map $ E: M_{\elli} \rightarrow M_{\elli} $. Moreover the relation $ E^2 = E $ holds 
by construction. The map  $ E $ is the natural projection onto the $ T $-invariant 
elements in $ M_{\elli} $. 

\section{Equivariant periodic cyclic homology}\label{secper}

In this section we recall the definition of bivariant equivariant periodic cyclic homology $ HP^G_*(A,B) $ 
given in \cite{Voigtperiodic}. 
An important property of the construction is that the fundamental object in the theory, the equivariant 
Hodge tower $ \theta \Omega_G(A) $ of a $ G $-algebra $ A $ is not a complex in the usual sense. 
In fact, the differential $ \partial $ in $ \theta \Omega_G(A) $ fails 
to satisfy $ \partial^2 = 0 $. \\
First we have to discuss noncommutative equivariant differential forms. 
Let $ A $ be a $ G $-algebra. The equivariant $ n $-forms of $ A $ are defined by 
$ \Omega^n_G(A) = \mathcal{O}_G \cotimes \Omega^n(A) $ where $ \Omega^n(A) = A^+ \cotimes A^{\cotimes n} $ and 
$ A^+ $ denotes the unitarization of $ A $. The group $ G $ acts diagonally on $ \Omega^n_G(A) $ and we have 
an obvious $ \mathcal{O}_G $-module structure. In this way $ \Omega^n_G(A) $ becomes a covariant module. \\
The equivariant Hochschild boundary $ b: \Omega^n_G(A) \rightarrow \Omega^{n - 1}_G(A) $ is defined by 
\begin{align*}
b(f(t) \otimes &x_0 dx_1 \cdots dx_n) = f(t) \otimes x_0 x_1 dx_2 \cdots dx_n \\
&+ \sum_{j = 1}^{n - 1} (-1)^j f(t) \otimes x_0 dx_1 \cdots d(x_j x_{j + 1}) \cdots dx_n \\
&+ (-1)^n f(t) \otimes (t^{-1} \cdot x_n)x_0 dx_1 \cdots dx_{n - 1}.
\end{align*}
Moreover we have the equivariant Connes operator $ B: \Omega^n_G(A) \rightarrow \Omega^{n + 1}_G(A) $ which is given by 
\begin{equation*}
B(f(t) \otimes x_0dx_1 \cdots dx_n) = \sum_{i = 0}^n (-1)^{ni}
f(t) \otimes t^{-1} \cdot(dx_{n + 1 - i} \cdots dx_n)dx_0 \cdots
dx_{n - i}.
\end{equation*}
It is straightforward to check that $ b $ and $ B $ are covariant maps. The natural symmetry operator $ T $ 
for covariant modules is of the form
\begin{equation*}
T(f(t) \otimes \omega) = f(t) \otimes t^{-1} \cdot \omega 
\end{equation*}
on $ \Omega^n_G(A) $. One easily obtains the relations $ b^2 = 0, B^2 = 0 $ and $ Bb + bB = \id - T $ for 
these operators. This shows that $ \Omega_G(A) $ is a paramixed complex in the following sense. 
\begin{definition} A paramixed complex $ M $ of covariant modules is a sequence of
covariant modules $ M_n $ together with differentials $ b $ of degree $ -1 $ and $ B $ of degree $ + 1 $ satisfying
$ b^2 = 0 $, $ B^2 = 0 $ and
\begin{equation*}
[b,B] = bB + Bb = \id - T.
\end{equation*}
\end{definition}
The most important examples of paramixed complexes are bounded below in the sense that $ M_n = 0 $ if $ n < N $ 
for some fixed $ N \in \mathbb{Z} $. In particular, the equivariant differential forms $ \Omega_G(A) $ of a $ G $-algebra $ A $
satisfy this condition for $ N = 0 $. \\
The Hodge filtration of a paramixed complex $ M $ of covariant modules is defined by 
$$
F^n M = bM_{n + 1} \oplus \bigoplus_{j > n} M_j. 
$$
Clearly $ F^n M $ is closed under the operators $ b $ and $ B $. We write 
$$
L^n M = F^{n - 1}M/ F^n M
$$
for the $ n $th layer of the Hodge filtration. If $ M $ is bounded below such that $ M_n = 0 $ for $ n < 0 $ we define 
the $ n $th level $ \theta^n M $ of the Hodge tower of $ M $ 
by 
$$
\theta^n M = \bigoplus_{j = 0}^{n - 1} M_j \oplus M_n/b M_{n + 1}.
$$
By definition, the Hodge tower of $ M $ is the projective system $ \theta M = (\theta^n M)_{n \in \mathbb{N}} $. \\
We remark that the image of the Hochschild operator $ b $ is not closed in general. In this case the spaces $ F^nM $ are 
possibly incomplete and $ L^n M $ and $ \theta^nM $ may fail to be separated. However, in the examples that we will 
study the image of $ b $ is always closed and hence these problems do not show up. \\
The spaces $ \theta^nM $ are equipped with the grading into even and odd forms and the differential $ \partial = B + b $. 
In this way the Hodge tower becomes a projective system of paracomplexes in the following sense. 
\begin{definition}
A paracomplex of covariant modules is a $ \mathbb{Z}_2 $-graded covariant module $ C $ with a boundary operator 
$ \partial: C \rightarrow C $ of degree one such that $ \partial^2 = \id - T $. 
\end{definition}
Chain maps of paracomplexes and homotopy equivalences are defined by the usual formulas. \\
We will work with the following definition of equivariant periodic cyclic homology which is 
equivalent to the one given in \cite{Voigtperiodic} using $ X $-complexes.
\begin{definition}
Let $ G $ be a totally disconnected group and let $ A $ and $ B $ be $ G $-algebras. The bivariant equivariant periodic
cyclic homology of $ A $ and $ B $ is
\begin{equation*}
HP^G_*(A,B) =
H_*(\SHom_G(\theta\Omega_G(A \hat{\otimes} \mathcal{K}_G), \theta \Omega_G(B \hat{\otimes} \mathcal{K}_G))).
\end{equation*}
\end{definition}
To explain this definition we first remark that the $ G $-algebra $ \mathcal{K}_G $ was defined in section \ref{sectotdis}. 
Secondly, the definition involves covariant maps between projective systems of covariant modules. 
Maps between projective systems are always understood in the sense 
$$
\Hom((M_i)_{i \in I}, (N_j)_{j \in J}) = \varprojlim_{j \in J} \varinjlim_{i \in I} \Hom(M_i,N_j)
$$
of pro-categories. Finally, we consider the usual differential for a $ \Hom $-complex given by 
$$
\partial(\phi) = \phi \partial_A - (-1)^{|\phi|} \partial_B \phi
$$
for a homogeneous element $ \phi $ in order to define homology. This makes sense since 
the failure of the individual differentials in $ \theta\Omega_G(A \hat{\otimes} \mathcal{K}_G) $ and 
$ \theta\Omega_G(B \hat{\otimes} \mathcal{K}_G) $ to satisfy $ \partial^2 = 0 $ is cancelled out by naturality of the 
operator $ T $. \\
According to the decomposition of a covariant module into an elliptic and a hyperbolic part discussed in section \ref{sectotdis} 
we obtain a decomposition 
$$ 
HP^G_*(A,B) = HP^G_*(A,B)_{\elli} \oplus HP^G_*(A,B)_{\hyp} 
$$ 
of equivariant periodic cyclic homology where by definition 
$$
HP^G_*(A,B)_{\elli} = H_*(\SHom_G(\theta\Omega_G(A \hat{\otimes} \mathcal{K}_G)_{\elli}, \theta \Omega_G(B \hat{\otimes} \mathcal{K}_G)_{\elli}))
$$
and accordingly for $ HP^G_*(A,B)_{\hyp} $. As we shall see, and as it is familiar from the computation of the cyclic homology of group algebras,
the elliptic and hyperbolic parts of equivariant cyclic homology behave differently. \\
Let us have a closer look at the elliptic part and recall the definition of the projection $ E $ from section \ref{sectotdis}. 
Since the map $ E: \theta \Omega_G(A)_{\elli} \rightarrow \theta \Omega_G(A)_{\elli} $ commutes 
with the boundary operators $ b $ and $ B $ we get a direct sum decomposition
\begin{equation*}
\theta\Omega_G(A)_{\elli} = E \theta\Omega_G(A)_{\elli} \oplus (1 - E) \theta\Omega_G(A)_{\elli}
\end{equation*}
of paracomplexes. Actually, $ E \theta\Omega_G(A)_{\elli} $ is an ordinary complex since $ T = \id $ on this space. 
\begin{prop} \label{Econtr}
Let $ A $ be a $ G $-algebra. The paracomplex $ (1 - E) \theta\Omega_G(A)_{\elli} $ is covariantly contractible. Hence 
the canonical projection $ \theta\Omega_G(A)_{\elli} \rightarrow E\theta\Omega_G(A)_{\elli} $ is a covariant homotopy equivalence. 
\end{prop}
\proof The covariant map $ \id - T: \theta\Omega_G(A)_{\elli} \rightarrow \theta\Omega_G(A)_{\elli} $ is diagonal with respect to the 
direct sum decomposition $ E \theta\Omega_G(A)_{\elli} \oplus (1 - E) \theta\Omega_G(A)_{\elli} $. It follows from the definitions 
that $ \id - T $ is injective on $ (1 - E)\theta\Omega_G(A)_{\elli} $ and $ \id - T = 0 $ on 
$ E \theta\Omega_G(A)_{\elli} $. Now consider an element $ (1 - E) m \in (1 - E) \theta^n\Omega_G(A)_{\elli} $ with 
$ m \in \theta^n\Omega_G(A)_{\elli} $. We find $ k \in \mathbb{N} $ such that $ T^k(m) = \id $ and 
hence 
$$
(1 - E)(m) = (\id - T) \sum_{j = 0}^{k - 1} \frac{k - j}{k} \; T^j(m)
$$
is in the image of $ \id - T $. It follows that the map $ \id - T: (1 - E) \theta\Omega_G(A)_{\elli} \rightarrow 
(1 - E) \theta\Omega_G(A)_{\elli} $ is an isomorphism. Since $ \id - T $ is homotopic to zero we 
conclude that $ (1 - E) \theta\Omega_G(A)_{\elli} $ is covariantly contractible. \qed 

\section{Smooth functions on simplicial complexes} \label{secsmooth}

In this section we study smooth functions and smooth differential forms on simplicial complexes. \\
First we have to fix some notation. We denote by $ \Delta^k $ the 
$ k $-dimensional standard simplex 
\begin{equation*}
\Delta^k = \{(x_0,\dots,x_k) \in \mathbb{R}^{k + 1}|\, 0 \leq x_j \leq 1, 
\sum_{j = 0}^k x_j = 1 \}
\end{equation*}
in $ \mathbb{R}^{k + 1} $. By construction $ \Delta^k $ is contained 
in a unique $ k $-dimensional affine subspace of $ \mathbb{R}^{k + 1} $ which will 
be denoted by $ A^k $. A function $ f: \Delta^k \rightarrow \mathbb{C} $ is called smooth if it is the restriction of a 
smooth function on the affine space $ A^k $. \\
To obtain an appropriate class of functions for our purposes we have to 
require conditions on the behaviour of such smooth functions near the boundary 
$ \partial \Delta^k $ of the simplex $ \Delta^k $. Roughly speaking, we shall 
consider only those functions which are constant in the direction orthogonal to 
the boundary in a neighborhood of $ \partial \Delta^k $. \\
Let us explain this precisely. We denote by $ \partial^i \Delta^k $ the $ i $-th face of the standard simplex 
consisting of all points $ (x_0,\cdots,x_k) \in \Delta^k $ satisfying 
$ x_i = 0 $. Then $ \partial^i \Delta^k $ defines a hyperplane 
$ A^k_i \subset A^k $ in a natural way. To this 
hyperplane we associate the vector space $ V_i $ which contains all 
vectors in $ \mathbb{R}^k $ that are orthogonal to $ A^k_i $.
For $ v \in \mathbb{R}^k $ denote by $ \partial_v(f) $ the partial derivative 
of a smooth function $ f $ on $ A^k $ in direction $ v $. 
We say that a smooth function $ f: \Delta^k \rightarrow \mathbb{C} $ is $ i $-regular if there exists a neighborhood $ U_i $ of 
$ \partial^i \Delta^k $ such that $ \partial_v(f)(x) = 0 $ for all 
$ x \in U_i $ and all $ v \in V_i $. If we want to emphasize the 
particular neighborhood $ U_i $ we also say that $ f $ is $ i $-regular on $ U_i $. 
The function $ f $ is called regular if it is $ i $-regular for all 
$ i = 0, \dots, k $. 
We denote by $ C^\infty(\Delta^k) $ the algebra of regular smooth functions on 
$ \Delta^k $. \\ 
The idea behind these definitions is as follows. Let us denote by 
$ C^\infty(\Delta^k,\partial \Delta^k) \subset C^\infty(\Delta^k) $ the subalgebra consisting of those functions that vanish 
on the boundary $ \partial \Delta^k $ of $ \Delta^k $. 
It is not hard to check that $ C^\infty(\Delta^k,\partial \Delta^k) $ can be 
identified with the algebra $ C^\infty_c(\Delta^k\setminus \partial\Delta^k) $ of smooth functions with 
compact support on the open set $ \Delta^k \setminus \partial \Delta^k $. 
Moreover the inclusion $ \partial^i: \partial^i\Delta^k \rightarrow \Delta^k $ of a face induces a 
homomorphism $ C^\infty(\Delta^k) \rightarrow 
C^\infty(\partial^i\Delta^k) $. Our definition yields a natural class of smooth functions satisfying 
these properties. \\
If we identify $ \Delta^1 $ with the unit interval $ [0,1] $ the algebra $ C^\infty(\Delta^1) $ corresponds 
to the algebra of smooth functions on $ [0,1] $ which are constant around the endpoints. 
Moreover $ C^\infty(\Delta^1,\partial\Delta^1) $ can be identified with
the algebra $ C^\infty_c(0,1) $ of smooth functions with compact support 
on the open interval $ (0,1) $. \\
We want to extend the definition of regular smooth functions to arbitrary simplicial complexes. 
A regular smooth function on a simplicial complex $ X $ is given 
by a family $ (f_\sigma)_{\sigma \subset X} $ of regular smooth functions on the simplices of $ X $ which 
is compatible with restriction to faces in the 
obvious way. The function $ f $ is said to have compact support if only finitely many 
$ f_\sigma $ in the corresponding family are different from zero. 
We denote by $ C^\infty_c(X) $ the algebra of regular smooth functions 
with compact support on $ X $. If the simplicial complex $ X $ is finite we simply 
write $ C^\infty(X) $ instead of $ C^\infty_c(X) $. \\
Let us now describe the natural locally convex topology 
on the algebra $ C^\infty_c(X) $ of regular smooth functions on the simplicial 
complex $ X $. First we consider again the case $ X = \Delta^k $. 
Let $ \mathcal{U} = (U_0,\cdots,U_k) $ be a family of open 
subsets of $ \Delta^k $ where each $ U_i $ is a neighborhood of 
$ \partial^i\Delta^k $. The collection of all such families is partially ordered 
where $ \mathcal{U} \prec \mathcal{V} $ iff 
$ V_j \subset U_j $ for all $ j $ in the corresponding families. 
For a family $ \mathcal{U} = (U_0,\cdots,U_k) $ we let 
$ C^\infty(\Delta^k, \mathcal{U}) \subset C^\infty(\Delta^k) $ be the 
subalgebra of smooth functions which are $ i $-regular on $ U_i $ for all $ i $. 
We equip $ C^\infty(\Delta^k, \mathcal{U}) $ with the natural 
Fr\'echet topology of uniform convergence of all derivatives on 
$ \Delta^k $. For $ \mathcal{U} \prec \mathcal{V} $ we have an obvious 
inclusion $ C^\infty(\Delta^k, \mathcal{U}) \subset C^\infty(\Delta^k, \mathcal{V}) $ 
which is compatible with the topologies. Moreover 
$ C^\infty(\Delta^k) $ is obtained as the union of the algebras 
$ C^\infty(\Delta^k, \mathcal{U}) $. We equip $ C^\infty(\Delta^k) $ 
with the resulting inductive limit topology.
In this way the algebra $ C^\infty(\Delta^k) $ becomes a nuclear LF-algebra 
and the natural restriction homomorphism $ C^\infty(\Delta^k) \rightarrow 
C^\infty(\partial^i\Delta^k) $ associated to the 
inclusion of a face is continuous. \\
In order to introduce a topology on $ C^\infty_c(X) $ for arbitrary $ X $ let 
$ K \subset X $ be a finite subcomplex. 
A function $ f = (f_\sigma)_{\sigma \subset X} \in C^\infty_c(X) $ is said to have support in $ K $ if $ f_\sigma = 0 $ for 
all simplices $ \sigma \subset X $ not 
contained in $ K $. We let $ C^\infty_K(X) \subset C^\infty_c(X) $ be the algebra 
of regular smooth functions with support in $ K $. 
The algebra $ C^\infty_K(X) $ is equipped with the subspace 
topology from the finite direct sum of algebras $ C^\infty(\sigma) $ for 
$ \sigma \subset K $. If $ K \subset L $ are finite subcomplexes 
the obvious inclusion $ C^\infty_K(X) \rightarrow C^\infty_L(X) $ is 
compatible with the topologies. Moreover $ C^\infty_c(X) $ is the union 
over all finite subcomplexes $ K $ of the algebras $ C^\infty_K(X) $. Hence we 
obtain a natural inductive limit topology on $ C^\infty_c(X) $. We equip 
the algebra $ C^\infty_c(X) $ with the associated precompact bornology 
(which is equal to the bounded bornology). In this way $ C^\infty_c(X) $ becomes a 
complete bornological algebra. \\
Let us have a closer look at the natural restriction homomorphism 
$ C^\infty(\Delta^k) \rightarrow C^\infty(\partial\Delta^k) $. 
As above we denote by  $ C^\infty(\Delta^k,\partial \Delta^k) \subset C^\infty(\Delta^k) $ 
the kernel of this homomorphism which consists of all regular smooth functions on 
$ \Delta^k $ that have compact support in the interior $ \Delta^k \setminus \partial\Delta^k $. 
\begin{prop}\label{SimplexSplitting} For all $ k $ the 
restriction homomorphism $ C^\infty(\Delta^k) \rightarrow 
C^\infty(\partial \Delta^k) $ has a continuous linear splitting. 
Hence we obtain a linearly split extension
\begin{equation*} 
    \xymatrix{
     C^\infty(\Delta^k, \partial \Delta^k) \;\; \ar@{>->}[r] &
        C^\infty(\Delta^k) \ar@{->>}[r] \ar@{<.}[l]& C^\infty(\partial \Delta^k)    
     }
\end{equation*}
of complete bornological algebras. 
\end{prop}
\proof By definition we have $ C^\infty(\partial\Delta^0) = 0 $ and hence the case 
$ k = 0 $ is trivial. \\
For $ k = 1 $ we identify $ \Delta^1 $ with the unit 
interval $ [0,1] $. Choose a smooth function $ h: [0,1] \rightarrow [0,1] $ 
such that $ h = 1 $ on $ [0,1/3] $ and $ h = 0 $ on $ [2/3,1] $. 
We set $ e_0 = h, e_1 = 1 - h $ and define $ \sigma_1: C^\infty(\partial\Delta^1) = 
\mathbb{C} \oplus \mathbb{C} \rightarrow C^\infty(\Delta^1) $ by 
$ \sigma_1(f_0,f_1) = f_0 e_0  + f_1 e_1 $. It is clear that $ \sigma_1 $ 
is a continuous linear splitting for the restriction map. \\
In order to treat the case $ k = 2 $ we first consider the 
corresponding lifting problem for a corner of $ \Delta^2 $ which can be formulated as follows. 
Let us write $ \mathbb{R}^+ $ for the set of nonnegative real numbers. 
A corner of $ \Delta^2 $ can be viewed as a neighborhood of 
the point $ (0,0) $ in $ \mathbb{R}^+ \times \mathbb{R}^+ $. 
Given smooth functions $ f_1, f_2: \mathbb{R}^+ \rightarrow \mathbb{C} $ that are both constant in a neighborhood of $ 0 $ 
and satisfy $ f_1(0) = f_2(0) $ we want to construct a smooth function
$ f: \mathbb{R}^+ \times \mathbb{R}^+ \rightarrow \mathbb{C} $ 
such that 
\begin{bnum}
\item[a)] $ f(x_1,0) = f_1(x_1) $ and $ f(0,x_2) = f_2(x_2) $ 
for all $ x_1,x_2 \in \mathbb{R}^+ $,
\item[b)] the function $ f $ is constant in the transversal direction in a neighborhood of the 
boundary $ (\mathbb{R}^+ \times \{0\}) \cup (\{0\} \times \mathbb{R}^+) $,
\item[c)] $ f $ depends linearly and continuously on $ f_1 $ and $ f_2 $.
\end{bnum}
In order to construct such a function we first extend 
$ f_1 $ and $ f_2 $ to smooth functions $ F_1 $ and $ F_2 $ on $ \mathbb{R}^+ \times \mathbb{R}^+ $ by setting 
\begin{equation*}
F_1(x_1,x_2) = f_1(x_1), \qquad F_2(x_1,x_2) = f_2(x_2). 
\end{equation*}
Then we use polar coordinates $ (r, \theta) $ in 
$ (\mathbb{R}^+ \times \mathbb{R}^+) \setminus \{ (0,0) \} $ to define a smooth function $ g_1 $ by 
\begin{equation*}
g_1(r, \theta) = h\biggl(\frac{2 \theta}{\pi}\biggr)
\end{equation*}
where $ h $ is the function from above. We extend $ g_1 $ to 
$ \mathbb{R}^+ \times \mathbb{R}^+ $ by setting $ g_1(0,0) = 0 $. 
Moreover we define $ g_2 $ by $ g_2 = 1 - g_1 $ on 
$ (\mathbb{R}^+ \times \mathbb{R}^+) \setminus\{(0,0)\} $ and $ g_2(0,0) = 0 $. 
Remark that $ g_1 $ and $ g_2 $ are 
not continuous in $ (0,0) $. We can now define the desired 
function $ f: \mathbb{R}^+ \times \mathbb{R}^+ \rightarrow \mathbb{C} $ by 
\begin{equation*}
f = f_1(0) \delta + F_1 g_1 + F_2 g_2 = f_2(0) \delta + F_1 g_1 + F_2 g_2
\end{equation*}
where $ \delta $ is the characteristic function of the point $ (0,0) $. 
The function $ f $ is smooth in $ (0,0) $ since the assumptions on $ f_1 $ and 
$ f_2 $ imply that $ f $ is constant in a neighborhood of this point. 
Moreover it is easy to verify that $ f $ satisfies the conditions a), b) and c) 
above. Hence this construction solves the lifting problem for a corner of 
$ \Delta^2 $. \\
Now we want to show that the restriction map $ C^\infty(\Delta^2) \rightarrow 
C^\infty(\partial \Delta^2) $ has a continuous linear splitting. 
One can combine the functions $ g_0 $ and $ g_1 $ constructed above 
for the corners of $ \Delta^2 $ to obtain functions $ e_j $ on $ \Delta^2 $ 
for $ j = 0,1,2 $ such that 
\begin{bnum}
\item[a)] each $ e_j $ is a regular smooth function on $ \Delta^2 $ except in the vertices where $ e_j $ is zero, 
\item[b)] $ e_j = 1 $ in the interior of $ \partial^j\Delta^2 $ and 
$ e_j = 0 $ in the interior of $ \partial^i \Delta^2 $ for $ i \neq j $, 
\item[c)] for each $ j $ there exists a neighborhood $ U_j $ of the $ j $-th 
vertex $ v_j $ of $ \Delta^2 $ such that 
\begin{equation*}
\sum_{i \neq j} e_i = 1 
\end{equation*}
holds on $ U_j $ except in $ v_j $. 
\end{bnum}
Now assume that a regular smooth function $ f = (f_0,f_1,f_2) $ on 
$ \partial \Delta^2 $ is given where $ f_i $ is defined on the face 
$ \partial^i\Delta^2 $. The functions $ f_i $ can be 
extended to $ i $-regular smooth functions $ F_i $ on $ \Delta^2 $ by 
\begin{align*}
F_0(x_0,x_1,x_2) &= f_0\biggl(x_1 + \frac{x_0}{2}, x_2 + \frac{x_0}{2}\biggr)\\ 
F_1(x_0,x_1,x_2) &= f_1\biggl(x_0 + \frac{x_1}{2}, x_2 + \frac{x_1}{2}\biggr)\\
F_2(x_0,x_1,x_2) &= f_2\biggl(x_0 + \frac{x_2}{2}, x_1 + \frac{x_2}{2}\biggr).
\end{align*}
Moreover let $ \chi: \Delta^2 \rightarrow \mathbb{C} $ be the characteristic 
function of the set $ \{v_0,v_1,v_2\}$ consisting of the three vertices of 
$ \Delta^2 $. 
Using these functions we define $ \sigma_2(f): \Delta^2 \rightarrow \mathbb{C} $ by 
\begin{equation*}
\sigma_2(f)= f \chi + F_0 e_0 + F_1 e_1 + F_2 e_2. 
\end{equation*}
To avoid confusion we point out that $ f \chi $ is the function 
which is equal to $ f $ in the vertices of $ \Delta^2 $ and extended by 
zero to the whole simplex. 
It is easy to see that the restriction of $ \sigma_2(f) $ to the boundary 
of $ \Delta^2 $ is equal to $ f $. Using the fact that 
$ F_i $ is $ i $-regular one checks that $ \sigma_2(f) $ is a regular smooth 
function on $ \Delta^2 $. Our construction yields a continuous linear map 
$ \sigma_2: C^\infty(\partial \Delta^2) \rightarrow C^\infty(\Delta^2) $ 
which splits the natural restriction homomorphism. \\
To prove the assertion for $ k > 2 $ one proceeds in a similar way as in the case 
$ k = 2 $. Essentially one has to combine the functions constructed above in an
appropriate way. 
First we consider again the lifting problem for a corner of 
$ \Delta^k $. Such a corner can be viewed as a neighborhood of $ (0,\dots,0) $ in 
$ (\mathbb{R}^+)^k $. 
We are given smooth functions $ f_1, \dots, f_k: (\mathbb{R}^+)^{k - 1} $ 
which are transversally constant in a neighborhood of the boundary and 
satisfy certain compatibility conditions. The function $ f_j $ is extended 
to a smooth function $ F_j: (\mathbb{R}^+)^k \rightarrow \mathbb{C} $ by setting
\begin{equation*}
F_j(x_1, \dots, x_k) = f_j(x_1,\dots,x_{j - 1},x_{j + 1}, \dots, x_k).
\end{equation*}
For $ 1 \leq i < j \leq k $ we define 
\begin{equation*}
g_{ij}(x_1,\dots,x_k) = g_1(x_i,x_j), \qquad g_{ji}(x_1,\dots,x_k) = g_2(x_i,x_j)
\end{equation*}
where $ g_1 $ and $ g_2 $ are the functions from above. 
Each function $ g_{ij} $ is smooth except in some 
$ k - 2 $-dimensional subspace inside the boundary and transversally 
constant in a neighborhood of the boundary. If we expand the 
product 
\begin{equation*}
\prod_{0 \leq i < j \leq k} (g_{ij} + g_{ji})
\end{equation*}
we obtain a sum of functions which are smooth except in the boundary and 
transversally constant in a neighborhood of the boundary. 
Moreover these functions vanish in the interior of all 
faces except possibly one. 
Using the functions $ F_j $ constructed before we can proceed as in the case 
$ k = 2 $ to solve the lifting problem for the $ k $-dimensional corner. 
Since this is a straightforward but lengthy verification we omit the details. \\
To treat the simplex $ \Delta^k $ we construct functions $ e_j $ for 
$ j = 0,\dots,k $ such that 
\begin{bnum}
\item[a)] $ e_j $ is regular smooth except in the $ (k - 2) $-skeleton of $ \Delta^k $ 
where $ e_j = 0 $,
\item[b)] $ e_j = 1 $ in the interior of $ \partial^j\Delta^k $ and $ e_j = 0 $ in 
the interior of $ \partial^i \Delta^k $ for $ i \neq j $, 
\item[c)] for each point $ x \in \partial \Delta^k $ there exists a neighborhood 
$ U_x $ of $ x $ such that 
\begin{equation*}
\sum_{j \in S(x)} e_j = 1
\end{equation*} 
in $ U_x $ except the $ (k - 2) $-skeleton of $ \Delta^k $ where 
$ S(x) $ is the collection of all $ i $ such that $ x \in \partial^i\Delta^k $. 
\end{bnum}
Assume that a regular smooth function $ f = (f_0,\dots ,f_k) $ on 
$ \partial \Delta^k $ is given where $ f_i $ is defined 
on the face $ \partial^i\Delta^k $. The function $ f_i $ can be 
extended to an $ i $-regular smooth function $ F_i $ on $ \Delta^k $ by 
\begin{equation*}
F_i(x_0,\dots, x_k) = f_i\biggl(x_0 + \frac{x_i}{k}, \dots, x_{i - 1}+ \frac{x_i}{k}, x_{i + 1} + \frac{x_i}{k}, \dots, x_k + \frac{x_i}{k} \biggr). 
\end{equation*}
Moreover let $ \chi: \Delta^k \rightarrow \mathbb{C} $ be the characteristic 
function of the $ k - 2 $-skeleton of $ \Delta^k $. We define 
$ \sigma_k(f): \Delta^k \rightarrow \mathbb{C} $ by 
\begin{equation*}
\sigma_k(f)= f \chi + \sum_{j = 0}^k F_j\, e_j.
\end{equation*}
Using the properties of the functions $ e_j $ and the fact that $ F_j $ is 
$ j $-regular one checks that $ \sigma_k(f) $ is a regular smooth function. 
The restriction of $ \sigma_k(f) $ to $ \partial\Delta^k $ is equal to $ f $. 
In this way we obtain a continuous linear 
splitting $ \sigma_k: C^\infty(\partial \Delta^k) \rightarrow C^\infty(\Delta^k) $ 
for the natural restriction homomorphism. \qed \\
For a simplicial complex $ X $ let $ X^k $ denote its $ k $-skeleton. 
Consider the natural continuous restriction homomorphism 
$ C^\infty_c(X^k) \rightarrow C^\infty_c(X^{k - 1}) $. 
The kernel of this homomorphism will be denoted by $ C^\infty_c(X^k,X^{k - 1}) $. 
\begin{prop}\label{ComplexSplitting} Let $ X $ be a simplicial complex. For all $ k $ 
the restriction homomorphism $ C^\infty_c(X^k) \rightarrow C^\infty_c(X^{k - 1}) $ 
has a continuous linear splitting. Hence we obtain a linearly split extension 
\begin{equation*} 
   \xymatrix{
      C^\infty(X^k,X^{k - 1}) \;\; \ar@{>->}[r] & C^\infty(X^k) \ar@{->>}[r] & C^\infty(X^{k - 1}) 
     }
\end{equation*}
of complete bornological algebras. 
\end{prop}
\proof We construct a retraction $ \rho: C^\infty_c(X^k) \rightarrow 
C^\infty_c(X^k,X^{k - 1}) $ for the natural inclusion. The algebra 
$ C^\infty_c(X^k,X^{k - 1}) $ can be identified with a direct sum 
$ \bigoplus_{i \in I} C^\infty(\Delta^k,\partial\Delta^k) $. 
Recall that the elements $ f \in C^\infty_c(X^k) $ are families 
$ (f_\sigma)_{\sigma \subset X^k} $. For each $ k $-simplex $ \eta \in X^k $ we 
define a map 
\begin{equation*}
\rho_\eta: C^\infty_c(X^k) \rightarrow C^\infty(\Delta^k,\partial\Delta^k),
\qquad \rho_\eta((f_\sigma)) = \rho_k(f_\eta)
\end{equation*}
where $ \rho_k: C^\infty(\Delta^k) \rightarrow 
C^\infty(\Delta^k,\partial \Delta^k) $ is the retraction obtained in proposition \ref{SimplexSplitting}. 
It is easy to check that $ \rho_\eta $ is continuous. 
The maps $ \rho_\sigma $ assemble to yield a map 
\begin{equation*}
\rho: C^\infty_c(X^k) \rightarrow \bigoplus_{i \in I}  C^\infty(\Delta^k,\partial\Delta^k) = C^\infty(X^k, X^{k - 1})
\end{equation*}
which is again continuous. Moreover by construction 
$ \rho $ is a retraction for the inclusion $ C^\infty(X^k, X^{k - 1}) 
\rightarrow C^\infty(X^k) $. \qed \\
We say that a complete bornological algebra $ K $ has local units if 
for every small subset $ S \subset K $ there exists an element $ e \in K $ 
such that $ e x = x e = x $ for all $ x \in S $. Clearly every 
unital complete bornological algebra has local units. 
In the bornological context the existence of local units has similar consequences as 
$ H $-unitality \cite{Wodzicki1}, \cite{Wodzicki2} in the algebraic setting. 
Clearly a complete bornological algebra which has local units is in particular 
$ H $-unital in the purely algebraic sense. 
The proof of excision in (algebraic) Hochschild homology for $ H $-unital algebras can 
immediately be adapted to show that every extension 
\begin{equation*} 
    \xymatrix{
     K \;\; \ar@{>->}[r] &
        E \ar@{->>}[r] \ar@{<.}[l]&
          Q   
     }
\end{equation*}
of complete bornological algebras with bounded linear splitting induces a long exact sequence in 
(bornological) Hochschild homology provided $ K $ has local units. A similar 
assertion holds for the homology with respect to the equivariant Hochschild boundary 
in the equivariant context.  
\begin{prop}\label{localunits} Let $ X $ be a locally finite simplicial complex. 
For every finite subcomplex $ K \subset X $ there exists a positive function 
$ e \in C^\infty_c(X) $ such that $ e = 1 $ on $ K $. In particular $ C^\infty_c(X) $ 
has local units. 
\end{prop}
\proof First recall that a simplicial complex $ X $ is called locally finite if every 
vertex of $ X $ is contained in only finitely many simplices of $ X $. 
A simplicial complex is locally finite iff it is a locally compact space 
in the weak topology. \\
The desired function $ e $ will be constructed inductively. On $ X^0 $ we define 
$ e(x) = 1 $ if $ x \in K^0 $ and $ e(x) = 0 $ otherwise. Assuming that $ e $ is constructed on $ X^{k - 1}$ we essentially 
have to extend functions which are defined on the boundary of $ k $-dimensional simplices to the whole simplices. 
If $ e $ is constant on the boundary we extend it to the whole simplex as a constant function. In general the extension can be done 
using the liftings for the restriction map $ C^\infty(\Delta^k) \rightarrow 
C^\infty(\partial \Delta^k) $ constructed in proposition \ref{SimplexSplitting}. 
It is clear that the resulting regular smooth function $ e $ is equal to $ 1 $ on 
$ K $. The fact that $ X $ is locally finite guarantees that $ e $ 
has compact support. Since every small subset of $ C^\infty_c(X) $ is contained in 
$ C^\infty_K(X) $ for some finite subcomplex $ K \subset X $ the previous 
discussion shows that $ C^\infty_c(X) $ has local units. \qed \\
Apart from smooth functions we also have to consider differential 
forms on simplicial complexes. 
A smooth differential form on the standard simplex $ \Delta^k $ is 
defined as the restriction of a smooth differential form on the $ k $-dimensional
affine space $ A^k $ to $ \Delta^k $. Again we have to impose some conditions on 
the behaviour near the boundary. Let us consider forms of a fixed degree $ p $. 
For $ v \in \mathbb{R}^k $ we denote by $ \mathcal{L}_v $ the Lie derivative 
in direction $ v $ and by $ \iota_v $ the interior product with the 
vector field associated to $ v $. 
Using the notation established in the beginning of this section we say that a 
smooth $ p $-form $ \omega $ on $ \Delta^k $ is $ i $-regular 
if there exist a neighborhood $ U_i $ of $ \partial^i\Delta^k $ such that 
$ \mathcal{L}_v(\omega)(x) = 0 $ and $ \iota_v(\omega)(x) = 0 $ 
for all $ x \in U_i $ and all $ v \in V_i $. The form $ \omega $ is called 
regular if it is $ i $-regular for all $ i = 0, \dots, k $. \\
Given a simplicial complex $ X $ a regular smooth $ p $-form 
$ \omega $ on $ X $ is a family $ (\omega_\sigma)_{\sigma \subset X} $ 
of regular smooth $ p $-forms on the simplices of $ X $ which is compatible with the natural restriction maps. 
A form $ \omega = (\omega_\sigma)_{\sigma \subset X} $ is said to have compact support if only finitely many 
$ \omega_\sigma $ in the corresponding family are nonzero. We denote by $ \mathcal{A}^p_c(X) $ the space 
of regular smooth $ p $-forms on $ X $ with compact support. 
The exterior differential $ d $ can be defined on $ \mathcal{A}_c(X) $ in the obvious 
way and turns it into a complex. Also the exterior product of 
differential forms extends naturally. 
Note that $ \mathcal{A}_c^0(X) $ can be identified with 
the algebra $ C^\infty_c(X) $ of regular smooth functions defined above. \\
As in the case of functions there is a natural topology on the space 
$ \mathcal{A}_c^p(X) $ of regular smooth $ p $-forms. Let us start with 
$ X = \Delta^k $ and consider a family $ \mathcal{U} = (U_0,\cdots,U_k) $ of open 
subsets of $ \Delta^k $ where each 
$ U_i $ is a neighborhood of $ \partial^i\Delta^k $. We let 
$ \mathcal{A}^p(\Delta^k, \mathcal{U}) \subset C^\infty(\Delta^k) $ be the 
space of smooth $ p $-forms which are $ i $-regular on $ U_i $ for all $ i $
and equip this space with the natural Fr\'echet topology. 
We obtain a corresponding inductive limit topology on $ \mathcal{A}^p(\Delta^k) $. 
Since one proceeds for an arbitrary simplicial complex $ X $ 
as in the case of functions we shall not work out the details. 
Most of the time we will not have to take into account the resulting bornology on
$ \mathcal{A}_c^p(X) $ in our considerations anyway. \\
We will have to consider differential forms not only as globally defined 
objects but also from the point of view of sheaf theory. 
The regularity conditions for smooth differential forms on a simplicial 
complex $ X $ obviously make sense also for an open subset $ U $ of $ X $. 
Hence we obtain in a natural way sheaves $ \mathcal{A}^p_X $ on $ X $ by 
letting $ \Gamma(U,\mathcal{A}^p_X) $ be the space of regular smooth 
$ p $-forms on the open set $ U \subset X $. We also write $ C^\infty_X $ 
for the sheaf $ \mathcal{A}^0_X $. The sheaf $ C^\infty_X $ is a sheaf of rings and the 
sheaves $ \mathcal{A}^p_X $ are sheaves of modules for $ C^\infty_X $. 
Clearly the space $ \Gamma_c(X,\mathcal{A}_X^p) $ of global sections with compact support of $ \mathcal{A}_X^p $ can be identified with 
$ \mathcal{A}_c^p(X) $. 
\begin{prop}\label{Resolutionprep} Let $ X $ be a locally finite simplicial complex. 
The sheaves $ \mathcal{A}^p_X $ are $ c $-soft for all $ p $ and 
\begin{equation*} 
    \xymatrix{
     \mathbb{C}_X \ar@{->}[r] &
        \mathcal{A}^0_X \ar@{->}[r]^d & \mathcal{A}^1_X \ar@{->}[r]^d & 
           \mathcal{A}^2_X \ar@{->}[r]^d & \cdots 
     }
\end{equation*}
is a resolution of the constant sheaf $ \mathbb{C}_X $ on $ X $. 
\end{prop}
\proof In this proof we will tacitly use some results from sheaf theory which can be found in \cite{Bredon2}. 
Let us first show that the sheaves $ \mathcal{A}^p_X $ are $ c $-soft. 
Since the sheaves $ \mathcal{A}^p_X $ are sheaves of modules for the sheaf of rings 
$ C^\infty_X $ it suffices to show that $ C^\infty_X $ 
is $ c $-soft. 
Using the fact that a sheaf $ \mathcal{F} $ on $ X $ is $ c $-soft iff the restrictions $ \mathcal{F}_{|K} $ are soft 
for all compact subsets $ K \subset X $ we may assume that $ X $ is a finite complex. 
We have to show that 
the restriction map $ \Gamma(X,C^\infty_X) \rightarrow \Gamma(K,C^\infty_X) $ 
is surjective for all closed subsets $ K \subset X $. Given a regular smooth 
function $ f $ on $ K $ we shall construct a regular smooth
function $ F: X \rightarrow \mathbb{C} $ which extends $ f $. 
For $ x \in X^0 $ we put $ F(x) = f(x) $ if $ x \in K $ and $ F(x) = 0 $ 
otherwise. Now assume that $ F $ has been constructed on $ X^{k - 1} $. 
In order to extend $ F $ to $ X^k $ we can consider each $ k $-simplex 
of $ X $ separately. If $ \sigma $ is a $ k $-simplex then 
$ F $ is already given on $ \partial \sigma $ by induction hypothesis and 
on the closed subset $ \sigma\cap K $ by assumption. The resulting function can 
be extended to a smooth regular function in a small neigborhood $ U $ of 
$ \partial \sigma \cup (\sigma\cap K) $. We find a regular smooth
function $ h $ on $ \sigma $ such that the support of $ h $ is contained in 
$ U $ and $ h = 1 $ on $ \partial \sigma \cup (\sigma\cap K) $. Using the function 
$ h $ we can extend $ F $ to the whole simplex $ \sigma $. \\
To show that the complex of sheaves 
$ \mathcal{A}^\bullet_X $ is a resolution of the constant sheaf on $ X $ we 
have to prove that the stalks $ (\mathcal{A}^\bullet_X)_x $ of this complex are  
resolutions of $ \mathbb{C} $ for all $ x \in X $. 
Each point $ x \in X $ is contained in 
$ X^k \setminus X^{k - 1} $ for some $ k $ and we find a 
$ k $-dimensional simplex $ \sigma $ in $ X $ such that $ x $ is an element 
in the interior $ \sigma \setminus \partial \sigma $ of $ \sigma $. 
From the definition of regular smooth differential forms we see 
that the stalks $ (\mathcal{A}^\bullet_X)_x $ depend only on the 
coordinates of $ \sigma $. Hence we can identify these stalks
in a natural way with stalks of the sheaves $ \mathcal{A}^\bullet_{\mathbb{R}^k} $ of smooth differential forms on $ k $-dimensional Euclidean 
space. The Poincar\'e lemma yields the assertion. \qed \\
For technical reasons we need a slightly more general class of spaces in the sequel. Namely, we will 
consider closed subspaces of spaces of the form $ T \times K $ where $ T $ is locally compact and totally disconnected 
and $ K $ is a locally finite simplicial complex. By definition, a regular smooth function $ f $ on the space $ T \times K $ 
is a function which can be written locally around any point as $ f(t,x) = F(x) $ for some 
regular smooth function $ F $ on $ K $. This condition clearly makes sense also for open subsets of $ T \times K $. 
Now let $ X \subset T \times K $ be a closed subspace. By definition, a regular smooth function on an open subset $ U \subset X $ is a 
function which can be extended to a regular smooth function in an open neighborhood of $ U $ in $ T \times K $. 
This defines a sheaf $ C^\infty_X $ of regular smooth functions on $ X $ in a natural way. Similarly one can consider 
differential forms and one obtains corresponding sheaves $ \mathcal{A}^p_X $. 
\begin{prop}\label{Resolution} Let $ T $ be a locally compact totally disconnected space 
and let $ K $ be a locally finite simplicial complex. For any closed subspace $ X \subset T \times K $ 
the sheaves $ \mathcal{A}^p_X $ are $ c $-soft and 
\begin{equation*} 
    \xymatrix{
     \mathbb{C}_X \ar@{->}[r] &
        \mathcal{A}^0_X \ar@{->}[r]^d & \mathcal{A}^1_X \ar@{->}[r]^d & 
           \mathcal{A}^2_X \ar@{->}[r]^d & \cdots 
     }
\end{equation*}
is a resolution of the constant sheaf $ \mathbb{C}_X $ on $ X $. 
\end{prop}
\proof The assertion concerning exactness can be proved as in proposition \ref{Resolutionprep}. 
It suffices to show that $ C^\infty_{T \times K} $ is $ c $-soft. Let $ A \subset T \times K $ be a compact subset 
and let $ f $ be a regular smooth function on $ A $. We have to construct a regular smooth function $ F $ on $ T \times K $ 
which extends $ f $. First choose open subsets $ U_1, \dots, U_n $ covering 
$ A $ with $ U_i = V_i \times W_i $ where $ V_i \subset T $ is compact open and $ W_i \subset K $ is open  
such that $ f $ can be extended to a regular smooth function 
on $ U = \bigcup_{i = 1}^n U_i $ which depends only on the $ K $-variable on each $ U_j $. 
In order to extend $ f $ to $ T \times K $ we proceed as follows. If we set $ V = \bigcup_{i = 1}^n V_i $ we find 
compact open subsets $ T_i \subset V $ such that the restriction of the extended function
$ f $ to $ T_i \times K \cap U $ does not 
depend on the $ T $-variable. Since the set $ T_j \times K \cap A $ is compact the same is true for its projection 
$ K_j = \pi_K(T_j \times K \cap A) $ to $ K $. By hypothesis $ g_j(x) = f(t,x) $ for $ (t,x) \in T_j \times K \cap A $ yields a well-defined 
regular smooth function on $ K_j $. We use proposition \ref{Resolutionprep} to extend $ g_j $ to a 
regular smooth function $ G_j $ on $ K $. Setting $ F(t,x) = G_j(x) $ for $ (t,x) \in T_j \times K $ we obtain a regular smooth function 
$ F $ on $ T_j \times K $ which restricts to 
$ f $ on $ T_j \times K \cap A $. Since $ T \setminus V \subset T $ is an open and closed subset 
we may set in addition $ F = 0 $ on $ (T \setminus V) \times K $ to obtain the desired extension of $ f $ to a regular smooth 
function on $ T \times K $. \qed 

\section{Bivariant equivariant cohomology}\label{secBSdef}

In this section we review the definition of bivariant equivariant cohomology given by Baum and Schneider \cite{BS}. \\
Let $ G $ be a totally disconnected group. A locally compact $ G $-space is a locally compact space $ X $ 
with a continuous action of $ G $. To every locally compact $ G $-space $ X $ we associate the Brylinski space 
\begin{equation*}
\hat{X} = \{(t,x) \in G \times X|\, t \,\;\text{is elliptic and}\,\; 
t \cdot x = x\} \subset G \times X
\end{equation*}
and the extended Brylinski space 
\begin{equation*}
\bar{X} = \{(t,x) \in G \times X|\; t \cdot x = x\} \subset G \times X.
\end{equation*}
Note that we have $ \hat{X} = \bar{X} $ if the action of $ G $ on $ X $ is proper.  
If $ G $ is discrete we may view $ \hat{X} $ as the disjoint union of 
the fixed point sets $ X^t = \{x \in X|\, t\cdot x = x\} $ of elements $ t \in G $ 
of finite order. Similarly, $ \bar{X} $ is the disjoint union of the fixed point 
sets $ X^t $ of arbitrary elements $ t \in G $ in this case. 
\begin{lemma} \label{closedlemma}
Let $ G $ be a totally disconnected group and let $ X $ be a locally compact $ G $-space. Then 
$ \hat{X} $ and $ \bar{X} $ are closed subspaces of $ G \times X $.   
\end{lemma}
\proof Let $ \mu: G \times X \rightarrow X \times X $ be the map defined by $ \mu(t,x) = (t\cdot x,x) $. 
Then $ \bar{X} = \mu^{-1}(\Delta) $ is the preimage of the diagonal 
$ \Delta \subset X \times X $ and hence closed. Since $ G_{\elli} $ is a closed subspace of $ G $ and 
$ \hat{X} = \bar{X} \cap (G_{\elli} \times X) $ it follows that $ \hat{X} $ is closed. \qed \\
There is a $ G $-action on the (extended) Brylinski space of a locally compact $ G $-space $ X $ given by the formula
\begin{equation*}
s \cdot (t,x) = (sts^{-1}, s\cdot x). 
\end{equation*}
In this way $ \hat{X} $ and $ \bar{X} $ become locally compact $ G $-spaces. \\
The space $ \bar{X} $ will appear in the equivariant Hochschild-Kostant-Rosenberg theorem 
in section \ref{secHKR}. For the remaining part of this section we will work only with 
the ordinary Brylinksi space $ \hat{X} $. \\
Since the category $ \Sh_G(\hat{X}) $ has enough injectives we can choose an injective resolution 
\begin{equation*} 
    \xymatrix{
     \mathbb{C}_{\hat{X}} \ar@{->}[r] &
        I^0 \ar@{->}[r] & I^1 \ar@{->}[r] & 
           I^2 \ar@{->}[r] & \cdots }
\end{equation*}
of the constant sheaf $ \mathbb{C}_{\hat{X}} $ in the category of equivariant sheaves on $ \hat{X} $. 
Consider the complex $ C_c^\bullet(\hat{X}) $ obtained by taking global sections with compact support of the sheaves $ I^\bullet $. 
Since the sheaves $ I^j $ are equivariant there is a natural $ G $-action 
on $ C_c^j(\hat{X}) $ for all $ j $. 
Moreover we have an $ \mathcal{O}_G $-module structure on $ C_c^j(\hat{X}) $ given by 
\begin{equation*}
(f \sigma)(s,x) = f(s) \sigma(s,x) 
\end{equation*}
for $ f \in \mathcal{O}_G $ and $ \sigma \in C_c^j(\hat{X}) $. Observe that only 
the elliptic part of $ \mathcal{O}_G $ acts nontrivially on $ C_c^j(\hat{X}) $. 
It is easy to check that this $ \mathcal{O}_G $-module structure and the natural 
$ G $-action combine to give each $ C_c^j(\hat{X}) $ the structure of a 
fine covariant module. \\
With these preparations the definition of bivariant equivariant cohomology 
given by Baum and Schneider can be formulated as follows. 
\begin{definition}\label{BSdef} Let $ G $ be a totally disconnected group and let $ X $ and $ Y $ be locally compact $ G $-spaces. 
The (delocalized) bivariant equivariant cohomology of $ X $ and $ Y $ is 
\begin{equation*}
H^n_G(X,Y) = \SExt_G^n(C^\bullet_c(\hat{X}),C^\bullet_c(\hat{Y}))
\end{equation*}
where $ \SExt_G $ denotes the hyperext functor in the category of fine covariant 
modules.  
\end{definition}
The functor $ \SExt_G $ can be viewed as the $ \Hom $-functor in the derived 
category of fine covariant modules. In order to compute the right-hand side in definition \ref{BSdef}
choose a complex $ I^\bullet(\hat{Y}) $ consisting of injective fine covariant modules 
together with a quasiisomorphism $ C^\bullet_c(\hat{Y}) \rightarrow 
I^\bullet(\hat{Y}) $. Then 
\begin{equation*}
\SExt_G^n(C^\bullet_c(\hat{X}), C^\bullet_c(\hat{Y})) = 
H_n(\SHom_G(C^\bullet_c(\hat{X}), I^\bullet(\hat{Y})),
\end{equation*}
hence in order to calculate $ \SExt_G $ we have to compute the homology of a 
certain $ \Hom $-complex. \\
Let us now specialize to group actions on simplicial complexes and give some 
more definitions. \\
Recall that a simplicial map between simplicial 
complexes $ X $ and $ Y $ is a continuous map $ f: X \rightarrow Y $ such that the restriction of 
$ f $ to any simplex of $ X $ is an affine map into a simplex of $ Y $. 
We say that the group $ G $ acts simplicially on $ X $ if every 
$ t \in G $ acts as a simplicial map. \\
Let $ G $ be a totally disconnected group. Assume that $ G $ acts simplicially on a simplicial complex $ X $. The 
action is called type-preserving if for each simplex $ \sigma $ of 
$ X $ the stabilizer $ G_\sigma $ fixes the vertices of $ \sigma $. In other words, 
an element of $ G $ which fixes a simplex actually acts trivially on this simplex. Passing to the barycentric subdivision one may always achieve
that $ G $ acts type-preserving. The action of $ G $ is called smooth if all isotropy groups are open. \\
Let us now specify the class of $ G $-spaces we are mainly interested in. 
\begin{definition} 
Let $ G $ be a totally disconnected group. A $ G $-simplicial complex is a 
simplicial complex $ X $ with a type-preserving smooth simplicial action of $ G $. 
A morphism of $ G $-simplicial complexes is an equivariant simplicial map $ f: X \rightarrow Y $. 
\end{definition}
Note that every $ G $-simplicial complex is a $ G $-$ CW $-complex. For the definition of 
a $ G $-$ CW $-complex we refer to \cite{LueckLNM}. If $ X $ is a $ G $-simplicial complex the space $ X^H $ of invariants with 
respect to a subgroup $ H \subset G $ is a subcomplex of $ X $. 
The action of $ G $ on $ X $ is proper iff the stabilizer of every point is a compact open subgroup of $ G $. \\
Let $ X $ and $ Y $ be locally finite $ G $-simplicial complexes. Our goal is to obtain a description of $ H^n_G(X,Y) $ 
which is closer to the definition of equivariant cyclic homology. \\
Due to lemma \ref{closedlemma} the Brylinski space $ \hat{X} $ is a closed subspace of $ G \times X $. Hence the formalism 
of regular smooth differential developped in section \ref{secsmooth} 
may be applied to $ \hat{X} $. According to proposition \ref{Resolution} 
we obtain a $ c $-soft resolution
\begin{equation*} 
    \xymatrix{
     \mathbb{C}_{\hat{X}} \ar@{->}[r] &
        \mathcal{A}^0_{\hat{X}} \ar@{->}[r]^d & \mathcal{A}^1_{\hat{X}} \ar@{->}[r]^d & 
           \mathcal{A}^2_{\hat{X}} \ar@{->}[r]^d & \cdots 
     }
\end{equation*}
of the constant sheaf $ \mathbb{C}_{\hat{X}} $ by regular smooth differential forms. This is a resolution in the category 
$ \Sh_G(\hat{X}) $ of equivariant sheaves on $ \hat{X} $. 
We equip the spaces $ \mathcal{A}_c^p(\hat{X}) $ of global sections with the fine bornology. 
\begin{prop} \label{BSC1} Let $ X $ and $ Y $ be locally finite $ G $-simplicial complexes. Then 
we have an isomorphism
\begin{equation*}
H^n_G(X,Y) = \SExt_G^n(\mathcal{A}^\bullet_c(\hat{X}), \mathcal{A}^\bullet_c(\hat{Y}))
\end{equation*}
which is natural with respect to equivariant proper simplicial maps in both variables. 
\end{prop}
\proof This isomorphism follows from proposition \ref{Resolution} and the fact that 
$ \SExt_G $ does not distinguish between quasiisomorphic complexes. 
The assertion concerning naturality is clear.  \qed \\
Now assume in addition that $ X $ is finite dimensional. Then the complex $ \mathcal{A}^\bullet_c(\hat{X}) $ is not only 
bounded below but also bounded above. This means that in order to compute 
$ \SExt_G^n(\mathcal{A}^\bullet_c(\hat{X}), \mathcal{A}^\bullet_c(\hat{Y})) $ we may use a complex 
$ P^\bullet(\hat{X}) $ consisting of projective fine covariant modules together 
with a quasiisomorphism 
$ p: P^\bullet(\hat{X}) \rightarrow \mathcal{A}^\bullet_c(\hat{X}) $
and obtain 
\begin{equation*}
\SExt_G^n(\mathcal{A}^\bullet_c(\hat{X}), \mathcal{A}^\bullet_c(\hat{Y})) = 
H^n(\SHom_G(P^\bullet(\hat{X}), \mathcal{A}^\bullet_c(\hat{Y}))).
\end{equation*}
If $ D = \dim(\hat{X}) $ is the dimension of $ \hat{X} $ we can construct 
a natural projective resolution $ P^\bullet(\hat{X}) $ of 
$ \mathcal{A}_c^\bullet(\hat{X}) $ in such a way that we obtain a commutative diagram 
of the form
\begin{equation*} 
    \xymatrix{
  \cdots \ar@{->}[r]^{\delta\quad\;} & P^{D - 2}(\hat{X}) \ar@{->}[r]^\delta \ar@{->}[d]^{p} & 
       P^{D - 1}(\hat{X}) \ar@{->}[r]^\delta \ar@{->}[d]^{p} &
       P^D(\hat{X}) \ar@{->}[r] \ar@{->}[d]^{p} & 0\\
  \cdots \ar@{->}[r]^{d\quad\;} & \mathcal{A}^{D - 2}_c(\hat{X}) \ar@{->}[r]^d & 
    \mathcal{A}^{D - 1}_c(\hat{X}) \ar@{->}[r]^d &
    \mathcal{A}^D_c(\hat{X}) \ar@{->}[r] & 0
     }
\end{equation*}
We may require in addition that $ P^j(\hat{X}) = EP^j(\hat{X}) $ for all $ j $ where 
$ E $ is the canonical projection on the elliptic part of a covariant module. This means in particular that the hyperbolic part of 
$ P^j(\hat{X}) $ is zero. In this case we call the projective resolution $ P^\bullet(\hat{X}) $ regular. 
Remark that $ T = \id $ for the natural operator $ T $ on a regular projective resolution $ P^\bullet(\hat{X}) $. \\
Let us view $ \mathcal{A}^\bullet_c(\hat{X}) $ as a mixed complex by setting the Hochschild boundary equal to zero 
and letting $ B = d $ be the exterior differential. To this mixed complex we associate a tower of supercomplexes 
$ \mathcal{A}_c(\hat{X}) = (\mathcal{A}_c(\hat{X})_k) $ as 
follows. We define
\begin{equation*}
\mathcal{A}_c(\hat{X})_k = 
\bigoplus_{j = 0}^k \mathcal{A}^j_c(\hat{X})  
\end{equation*}
and equip this space with the ordinary grading into even and odd forms 
and differential $ B + b = d $. Observe that $ \mathcal{A}_c(\hat{X})_k = \mathcal{A}_c(\hat{X})_D $ 
for $ k \geq D = \dim(\hat{X}) $. Hence the tower of supercomplexes $ \mathcal{A}_c(\hat{X}) $ 
is isomorphic to the constant supercomplex 
\begin{equation*}
\mathcal{A}_c(\hat{X}) \cong
\bigoplus_{j = 0}^D \mathcal{A}^j_c(\hat{X}).
\end{equation*}
In a similar way a regular projective resolution $ P^\bullet(\hat{X}) $ of $ \mathcal{A}^\bullet_c(\hat{X}) $ 
satisfies the axioms of a mixed complex. Let us define a tower of supercomplexes 
$ P(\hat{X}) = (P(\hat{X})_k) $ as follows.  
We set 
\begin{equation*}
P(\hat{X})_k = P^{-(k + 1)}(\hat{X})/\delta(P^{-(k + 2)}(\hat{X})) \oplus 
\bigoplus_{j = -k}^{k} P^j(\hat{X}) \oplus \delta(P^k(\hat{X}))
\end{equation*}
Remark that for $ k \geq D $ this becomes 
\begin{equation*}
P(\hat{X})_k = P^{-(k + 1)}(\hat{X})/\delta(P^{-(k + 2)}(\hat{X})) \oplus 
\bigoplus_{j = -k}^{D} P^j(\hat{X}).
\end{equation*}
Clearly we consider the grading into even and odd components on 
$ P(\hat{X})_k $ and equip these spaces with the differential $ \delta $. 
Recall from \cite{Voigtperiodic} that a a covariant pro-module is called relatively projective if it has the 
lifting property with respect to covariant maps between covariant pro-modules having a pro-linear section. 
Since the covariant modules $ P^j(\hat{X}) $ are projective for all $ j $ it 
is easy to see that the inverse system $ P(\hat{X}) $ is relatively projective. 
The chain map $ p: P^\bullet(\hat{X}) \rightarrow \mathcal{A}^\bullet_c(\hat{X}) $ 
induces a covariant chain map of supercomplexes $ p: P(\hat{X}) \rightarrow \mathcal{A}_c(\hat{X}) $. 
\begin{prop}\label{BSCProp1} Let $ X $ and $ Y $ be finite dimensional locally finite $ G $-simplicial complexes. Then
\begin{align*}
\bigoplus_{j \in \mathbb{Z}} H^{* + 2j}_G&(X,Y) 
= H_*(\varinjlim_k 
\SHom_G(P(\hat{X})_k, \mathcal{A}_c(\hat{Y}))) \\
&= H_*(\SHom_G(P(\hat{X}),\mathcal{A}_c(\hat{Y})))
\end{align*}
where in the last expression we take homomorphisms in the pro-category $ \pro(G\LSMod) $ of covariant 
modules. 
\end{prop}
\proof The component of degree $ n $ in  
$ \SHom_G(P^\bullet(\hat{X}), \mathcal{A}^\bullet_c(\hat{Y})) $ is 
\begin{equation*}
\bigoplus_{i \in \mathbb{Z}} \SHom_G(P^i(\hat{X}), \mathcal{A}^{i + n}_c(\hat{Y})),
\end{equation*}
here a direct sum occurs because $ \mathcal{A}^\bullet_c(\hat{Y}) $ is a bounded 
complex. We deduce 
\begin{equation*}
\bigoplus_{j \in \mathbb{Z}} 
\SHom_G^{* + 2j}(P^\bullet(\hat{X}), \mathcal{A}^\bullet_c(\hat{Y})) = 
\bigoplus_{j \in \mathbb{Z}} \bigoplus_{i \in \mathbb{Z}} \SHom_G(P^i(\hat{X}), \mathcal{A}^{i + 2j + *}_c(\hat{Y}))
\end{equation*}
and obtain natural maps 
\begin{equation*}
\lambda_k: \SHom_G(P(\hat{X})_k, \mathcal{A}_c(\hat{Y})) 
\rightarrow \bigoplus_{j \in \mathbb{Z}} 
\SHom_G^{* + 2j}(P^\bullet(\hat{X}), \mathcal{A}^\bullet_c(\hat{Y}))
\end{equation*}
for all $ k \geq \dim(\hat{X}) $. It is easy to check that each $ \lambda_k $ is a chain map. 
Moreover the maps $ \lambda_k $ are compatible 
with the projections in the first variable. 
The resulting map 
\begin{equation*}
\lambda: \varinjlim_k\SHom_G(P(\hat{X})_k, \mathcal{A}_c(\hat{Y}))
\rightarrow \bigoplus_{j \in \mathbb{Z}} 
\SHom_G^{* + 2j}(P^\bullet(\hat{X}), \mathcal{A}^\bullet_c(\hat{Y}))
\end{equation*}
is an isomorphism of complexes. \qed 

\section{The equivariant Hochschild-Kostant-Rosenberg theorem} \label{secHKR}

The algebra $ C^\infty_c(X) $ of regular smooth functions on a $ G $-simplicial complex $ X $ is a $ G $-algebra in a natural 
way. In this section we identify the homology of 
$ \Omega_G(C^\infty_c(X)) $ with respect to the equivariant Hochschild boundary. 
This will be an important ingredient in the proof of theorem \ref{BSC} below. \\
Recall from section \ref{secBSdef} the definition of the extended Brylinski space $ \bar{X} $. 
Let us view $ \mathcal{A}_c(\bar{X}) $ as a (para-) mixed complex with 
$ b $-boundary equal to zero and $ B $-boundary equal to the exterior 
differential $ d $. We define the equivariant Hochschild-Kostant-Rosenberg map 
$ \alpha: \Omega_G(C^\infty_c(X)) \rightarrow \mathcal{A}_c(\bar{X}) $ by 
\begin{equation*}
\alpha(f(t) \otimes a_0 da_1 \cdots da_n) = \frac{1}{n!}\; f(t)\; a_0 da_1 \wedge \cdots \wedge da_{n|X^t}
\end{equation*}
where we recall that $ X^t $ denotes the set of fixed points under the action of 
$ t $.
\begin{theorem} \label{HKR} Let $ G $ be a totally disconnected group and let $ X $ be a $ G $-simplicial complex. 
The equivariant Hochschild-Kostant-Rosenberg map 
$$
\alpha: \Omega_G(C^\infty_c(X)) \rightarrow \mathcal{A}_c(\bar{X}) 
$$ 
is a map of paramixed complexes and induces an isomorphism on the homology with respect to the Hochschild boundary. 
\end{theorem}
\proof Let us first show that $ \alpha $ is a map of paramixed complexes. 
We compute 
\begin{align*}
\alpha b&(f(t) \otimes a_0 da_1 \cdots da_n) = 
\sum_{j = 0}^{n - 1} (-1)^j \alpha(f(t) \otimes a_0 da_1 \cdots d(a_j a_{j + 1}) \cdots da_n) \\
&\qquad + (-1)^n \alpha(f(t) \otimes (t^{-1} \cdot a_n) a_0 da_1 \cdots 
da_{n - 1}) \\
&= \frac{1}{(n - 1)!}\;  \biggl(\,\sum_{j = 0}^{n - 1} (-1)^j 
f(t)\; a_0 da_1 \wedge \cdots \wedge d(a_j a_{j + 1}) \cdots \wedge da_{n|X^t} \\
&\qquad + (-1)^n f(t) a_n a_0 da_1 \wedge \cdots \wedge da_{n - 1|X^t}\biggr) = 0
\end{align*}
where we use $ (t^{-1} \cdot a_n)(x) = a_n(t \cdot x) = a_n(x) $ for all 
$ x \in X^t $. Moreover we have
\begin{align*}
\alpha B&(f(t) \otimes a_0 da_1 \cdots da_n) = 
\sum_{j = 0}^n (-1)^{nj} \alpha(f(t) \otimes t^{-1} \cdot (da_{n - j + 1} \cdots da_n) da_0 \cdots da_{n - j}) \\
&= \frac{1}{(n + 1)!}\; \sum_{j = 0}^n (-1)^{nj} f(t) \,
da_{n - j + 1} \wedge \cdots \wedge da_n \wedge da_0 \wedge \cdots 
\wedge da_{n - j| X^t} \\
&= \frac{1}{n!}\; f(t)\,
da_0 \wedge \cdots \wedge da_{n|X^t} = d \alpha(f(t) \otimes a_0 da_1 \cdots da_n)
\end{align*}
and hence $ \alpha $ commutes with the boundary operators as claimed. \\
In order to show that $ \alpha $ induces an isomorphism in homology it suffices to prove that the localized maps 
$$
\alpha_t: \Omega_G(C^\infty_c(X))_t \rightarrow \mathcal{A}_c(\bar{X})_t = \mathcal{A}_c(X^t)
$$
are quasiisomorphisms for all $ t \in G $. Let us consider 
the case that $ X $ is an equivariant simplex. By definition, an equivariant simplex is a space of the form 
$ X = G/H \times \Delta^k $  where $ H $ is an open subgroup of $ G $ and the action on 
$ G/H $ is given by translation. The boundary $ \partial X $ of $ X $ is defined 
by $ \partial X = G/H \times \partial \Delta^k $. Recall from section \ref{secsmooth} that $ C^\infty_c(X,\partial X) $ denotes the kernel 
the restriction map $ C^\infty_c(X) \rightarrow C^\infty(\partial X) $. 
Similarly, $ \mathcal{A}_c(X^t, \partial X^t) $ is the kernel of the natural map 
$ \mathcal{A}_c(X^t) \rightarrow \mathcal{A}_c(\partial X^t) $. The localized 
equivariant Hochschild-Kostant-Rosenberg map restricts to a chain map 
\begin{equation*}
\alpha_t: \Omega_G(C^\infty_c(X,\partial X))_t \rightarrow \mathcal{A}_c(X^t, \partial X^t).
\end{equation*}
Let us specialize further to the case $ G = \mathbb{Z} $ and $ t = 1 $. We write 
$ \Delta^k[n] $ for the $ \mathbb{Z} $-equivariant simplex $ \mathbb{Z}/n \mathbb{Z} \times \Delta^k $. Remark that the fixed 
point set $ \Delta^k[n]^1 $ for the action of $ 1 \in \mathbb{Z} $ is empty for $ n = 0 $ or $ n > 1 $.  
\begin{prop} \label{HKROrbitslokal} With the notation as above we have for every $ k \geq 0 $: 
\begin{bnum}
\item[a)] The localized equivariant Hochschild-Kostant-Rosenberg map 
\begin{equation*}
\alpha_1: \Omega_\mathbb{Z}(C^\infty(\Delta^k[1], \partial \Delta^k[1]))_1 \rightarrow \mathcal{A}(\Delta^k[1]^1, \partial \Delta^k[1]^1) 
\end{equation*}
is a quasiisomorphism. 
\item[b)] For $ n = 0 $ and $ n > 1 $ the homology of 
$ \Omega_\mathbb{Z}(C^\infty_c(\Delta^k[n], \partial \Delta^k[n]))_1 $ with respect to the 
equivariant Hochschild boundary is trivial. 
\end{bnum}
\end{prop}
\proof $ a) $ By definition we have $ \Delta^k[1] = \Delta^k $ and the action is trivial. 
The algebra $ C^\infty(\Delta^k[1], \partial \Delta^k[1]) $ can be identified with the 
algebra $ C^\infty_c(\Delta^k \setminus \partial \Delta^k) $ of smooth functions with compact support on $ \Delta^k\setminus\partial\Delta^k $. 
Moreover the space $ \mathcal{A}(\Delta^k[1]^1, \partial \Delta^k[1]^1) $ 
consists of differential forms with compact support on 
$ \Delta^k\setminus\partial\Delta^k $. Hence the assertion follows from the ordinary Hochschild-Kostant-Rosenberg theorem \cite{Connes1}, 
\cite{Teleman}. \\
$ b) $ Let $ n = 0 $ or $ n > 1 $ and let $ B $ be any unital complete bornological algebra. We equip 
$ B $ with the trivial $ \mathbb{Z} $-action and consider the $ \mathbb{Z} $-algebra 
$ C_c(\mathbb{Z}/n \mathbb{Z}) \hat{\otimes} B $. An element of this algebra can 
be written as a linear combination of elements $ x[i] $ 
where $ x[i] \in C_c(\mathbb{Z}/n\mathbb{Z}) \cotimes B  $ for $ i \in \mathbb{Z}/n\mathbb{Z} $ denotes the 
characteristic function located in $ i $ with value $ x \in B $. \\
We view the localized Hochschild complex $ \Omega_\mathbb{Z}(C_c(\mathbb{Z}/n \mathbb{Z}) \hat{\otimes} B)_1 $ as a double complex with two columns. 
This corresponds to the natural decomposition 
$$ 
\Omega^n_G(C) = \mathcal{O}_G \hat{\otimes} C^{\hat{\otimes} n + 1} \oplus \mathcal{O}_G \hat{\otimes} C^{\hat{\otimes} n} 
$$ 
of the space of equivariant differential forms of a $ G $-algebra $ C $. Since $ B $ is unital the algebra 
$ C_c(\mathbb{Z}/n \mathbb{Z}) \hat{\otimes} B $ has local units. Consequently, the natural inclusion 
of the first column of $ \Omega_\mathbb{Z}(C_c(\mathbb{Z}/n \mathbb{Z}) \hat{\otimes} B)_1 $ into the total complex induces an isomorphism in homology. 
Let us construct a contracting homotopy $ h $ for the first column of $ \Omega_\mathbb{Z}(C_c(\mathbb{Z}/n \mathbb{Z}) \hat{\otimes} B)_1 $ as follows.
We define 
\begin{align*}
h(x_0[i_0]& dx_1[i_1]\cdots dx_p[i_p]) = (-1)^{l + 1} x_0[i_0] dx_1[i_1] \cdots d1[i_{l + 1}] dx_{l + 1}[i_{l + 1}] \cdots dx_p[i_p]
\end{align*}
if $ 0 \leq l \leq p - 1 $ is the smallest number such that 
$ i_l \neq i_{l + 1} $ and
\begin{equation*}
h(x_0[i_0] dx_1[i_1] \cdots dx_p[i_p]) = 1[i_0] dx_0[i_0] \cdots dx_p[i_p] 
\end{equation*}
if $ i_j = i_{j + 1} $ for $ 0 \leq j \leq p - 1 $. An easy calculation shows that $ h $ is indeed a contracting homotopy. 
Hence the complexes $ \Omega_\mathbb{Z}(C_c(\mathbb{Z}/n \mathbb{Z}) \hat{\otimes} B)_1 $ are acyclic 
for $ n = 0 $ or $ n > 1 $. \\
Due to proposition \ref{ComplexSplitting} we have an extension of $ \mathbb{Z} $-algebras with bounded linear splitting
\begin{equation*} 
 \xymatrix{
     C^\infty_c(\Delta^k[n], \partial\Delta^k[n]) \;\; \ar@{>->}[r] &
        C^\infty_c(\Delta^k[n]) \ar@{->>}[r] \ar@{<.}[l]&  C^\infty_c(\partial\Delta^k[n])   
     }
\end{equation*}
The algebras $ C^\infty_c(\Delta^k[n]) $ and $  C^\infty_c(\partial\Delta^k[n]) $ 
are of the form $ C_c(\mathbb{Z}/n \mathbb{Z}) \hat{\otimes} B $ described above. Hence the complexes 
$ \Omega_\mathbb{Z}( C^\infty_c(\Delta^k[n]))_1 $ 
and $ \Omega_\mathbb{Z}( C^\infty_c(\partial\Delta^k[n]))_1 $ are acyclic. 
Since the algebra $ C^\infty_c(\Delta^k[n],\partial \Delta^k[n]) $ has local units 
we obtain a long exact sequence in homology showing that  
$ \Omega_\mathbb{Z}(C^\infty_c(\Delta^k[n], \partial \Delta^k[n])_1 $ is acyclic as well. This 
yields the claim. \qed \\ 
Let us come back to the localized Hochschild-Kostant-Rosenberg map for 
arbitrary $ G $ and $ t $ and an equivariant simplex $ X = G/H \times \Delta^k $. 
We extend proposition \ref{HKROrbitslokal} to this situation as follows. There is a canonical group homomorphism 
$ \mathbb{Z} \rightarrow G $ which maps $ 1 $ to $ t $ and we may view $ X = G/H \times \Delta^k $ as a 
$ \mathbb{Z} $-space in this way. 
Clearly the localized complexes $ \Omega_G(C^\infty_c(X,\partial X))_t $ and 
$ \Omega_\mathbb{Z}(C^\infty_c(X, \partial X))_1 $ are isomorphic since 
the equivariant Hochschild boundary in $ \Omega_G(C^\infty_c(X,\partial X))_t $ only depends on the action of $ t $. 
Viewed as a $ \mathbb{Z} $-space, $ X $ can be written as disjoint union 
\begin{equation*}
X = \bigcup_{j \in J} \Delta^k[n_j] 
\end{equation*}
for some index set $ J $ where $ \Delta^k[n] = \mathbb{Z}/n\mathbb{Z} \times \Delta^k $ as before. 
In this decomposition the spaces $ \Delta^k[n] $ may appear with multiplicity. Let us determine 
how $ \alpha_t: \Omega_\mathbb{Z}(C^\infty_c(X,\partial X))_t \rightarrow 
\mathcal{A}_c(X^t, \partial X^t) $ can be described in terms of the spaces 
$ \Delta^k[n_j] $. 
On the right hand side the decomposition of $ X $ induces a direct sum 
decomposition 
$$ 
\mathcal{A}_c(X^t, \partial X^t) = \bigoplus_{j \in J} \mathcal{A}_c(\Delta^k[n_j]^1, \partial\Delta^k[n_j]^1).
$$
Moreover we have an isomorphism 
$$
C^\infty_c(X) = \bigoplus_{j \in J} C^\infty_c(\Delta^k[n_j]) 
$$ 
of $ \mathbb{Z} $-algebras and a natural inclusion of complexes 
\begin{equation*} 
\iota: \bigoplus_{j \in J} 
\Omega_{\mathbb{Z}}(C^\infty_c(\Delta^k[n_j],\partial \Delta^k[n_j]))_1 \rightarrow 
\Omega_{\mathbb{Z}}\Bigl(\bigoplus_{j \in J} C^\infty_c(\Delta^k[n_j], \partial \Delta^k[n_j])\Bigr)_1
\end{equation*}
on the left hand side. It follows from the existence of local units and an inductive limit argument 
that the map $ \iota $ is a quasiisomorphism with respect to the Hochschild boundary. \\
Hence, up to quasiisomorphism, the map $ \alpha_t $ can be decomposed as a direct sum of maps 
$ \alpha_1: \Omega_\mathbb{Z}(C^\infty_c(\Delta^k[n_j], \partial \Delta^k[n_j]))_1 \rightarrow \mathcal{A}_c(\Delta^k[n_j]^1, \partial\Delta^k[n_j]^1) $.
We apply proposition \ref{HKROrbitslokal} and obtain the following statement. 
\begin{prop} \label{HKROrbits} 
Let $ G $ be a totally disconnected group and let $ X = G/H \times \Delta^k $ be an equivariant simplex. 
For every $ t \in G $ the localized equivariant Hochschild-Kostant-Rosenberg map 
\begin{equation*}
\alpha_t: \Omega_G(C^\infty_c(X,\partial X))_t \rightarrow \mathcal{A}_c(X^t,\partial X^t) 
\end{equation*}
is a quasiisomorphism.
\end{prop}
Let us now finish the proof of theorem \ref{HKR}. First we assume that the $ G $-simplicial complex $ X $ is finite dimensional and use 
induction on the dimension of $ X $. 
If $ \dim(X) = 0 $ the space $ X $ is a disjoint union of homogenous spaces $ G/H $. As above it suffices to consider a single equivariant simplex 
$ X = G/H $. Since we have $ C^\infty_c(X) = C^\infty_c(X,\partial X) $ in this case 
the assertion follows from proposition \ref{HKROrbits}. Assume that $ \dim(X) = k $ and that the assertion is proved for all 
$ G $-simplicial complexes of dimension $ k - 1 $. We consider the commutative diagram 
\begin{equation*} 
    \xymatrix{
    \Omega_G(C^\infty_c(X,X^{k - 1}))_t \ar@{->}[r] \ar@{->}[d]^{\alpha_t} & 
      \Omega_G(C^\infty_c(X))_t \ar@{->}[r] \ar@{->}[d]^{\alpha_t} &
       \Omega_G(C^\infty_c(X^{k - 1}))_t \ar@{->}[d] ^{\alpha_t}\\
 \mathcal{A}_c(X^t,(X^{k - 1})^t) \ar@{->}[r] & 
     \mathcal{A}_c(X^t) \ar@{->}[r] &
  \mathcal{A}_c((X^{k - 1})^t)  
     }
\end{equation*}
where $ X^{k - 1} $ denotes the $ k - 1 $-skeleton of $ X $. 
The algebra $ C^\infty_c(X, X^{k - 1}) $ is a direct sum of 
algebras of the form $ C^\infty_c(\sigma, \partial\sigma) $ where 
$ \sigma = G/H \times \Delta^k $ is an equivariant simplex. In particular
$ C^\infty_c(X^k, X^{k - 1}) $ has local units. Hence the upper horizontal sequence induces a long exact sequence 
in homology. Proposition \ref{HKROrbits} implies that the 
left vertical map is a quasiisomorphism. The right vertical map 
is a quasiisomorphism by induction hypothesis. Hence $ \alpha_t: \Omega_G(C^\infty_c(X))_t
\rightarrow \mathcal{A}_c(X^t) $ is a quasiisomorphism as well. \\
For an arbitrary $ G $-simplicial complex $ X $ we take the inductive limit over all finite dimensional subcomplexes in order to 
obtain the assertion. This completes the proof of theorem \ref{HKR}. 

\section{The comparison theorem} \label{secBSC}

In this section we prove the following theorem which describes the relation between equivariant periodic cyclic homology and 
bivariant equivariant cohomology in the sense of Baum and Schneider. 
\begin{theorem}\label{BSC} Let $ G $ be a totally disconnected group and let $ X $ and $ Y $ be finite dimensional locally 
finite $ G $-simplicial complexes. If the action of $ G $ on $ X $ is proper there exists a 
natural isomorphism
\begin{equation*}
HP^G_*(C^\infty_c(X), C^\infty_c(Y)) \cong 
\bigoplus_{j \in \mathbb{Z}} H^{* + 2j}_G(X,Y).
\end{equation*}
Even if the action on $ X $ is not proper the elliptic part of 
$ HP^G_*(C^\infty_c(X), C^\infty_c(Y)) $ is naturally isomorphic to $ \bigoplus_{j \in \mathbb{Z}} H^{* + 2j}_G(X,Y) $. 
These isomorphisms are natural with respect to equivariant proper simplicial 
maps in both variables. 
\end{theorem}
It follows immediately from the definitions that the theory defined by Baum and Schneider a priori 
only has an elliptic part. Hence theorem \ref{BSC} states in particular that the hyperbolic 
part of $ HP^G_*(C^\infty_c(X), C^\infty_c(Y)) $ is zero provided the action of 
$ G $ on $ X $ is proper. In general the hyperbolic part 
of $ HP^G_*(C^\infty_c(X), C^\infty_c(Y)) $ might be 
different from zero, however, this cannot be detected using the theory of Baum and Schneider. \\
The proof of theorem \ref{BSC} is divided into several steps. First we shall identify $ HP^G_* $ with an auxiliary bivariant theory 
$ h_*^G $ under the assumptions of the theorem. 
We denote by $ \mathfrak{Fine} $ the natural forgetful functor 
on covariant modules which changes the bornology to the fine bornology. 
The functor $ \mathfrak{Fine} $ is extended to the category $ \pro(G\LSMod) $ in 
the obvious way. With this notation we define the bivariant theory $ h^G_*(A,B) $ for 
$ G $-algebras $ A $ and $ B $ by 
\begin{equation*}
h^G_*(A,B) = 
H_*(\SHom_G(\mathfrak{Fine}(\theta \Omega_G(A \hat{\otimes} \mathcal{K}_G)), 
\mathfrak{Fine}(\theta \Omega_G(B \hat{\otimes} \mathcal{K}_G)))).
\end{equation*}
This definition is identical to the definition of $ HP^G_* $ except that we do 
not require the covariant maps in the $ \Hom $-complex to be bounded. 
Evidently $ h^G_* $ shares many properties with $ HP^G_* $. For our purposes it is important that $ h^G_*$ satisfies excision 
in both variables. This follows immediately from the generalized excision theorem in equivariant periodic cyclic 
homology \cite{Voigtperiodic}. Moreover there is an obvious composition product 
for $ h^G_* $ and a natural transformation 
\begin{equation*}
v: HP^G_*(A,B) \rightarrow h^G_*(A,B)
\end{equation*}
which is obtained by forgetting the bornology. It is clear that $ v $ is compatible with the composition product. 
\begin{prop} \label{hprop} Let $ X $ be a finite dimensional $ G $-simplicial complex and let $ B $ be an arbitrary 
$ G $-algebra. Then the natural map 
\begin{equation*}
v: HP^G_*(C^\infty_c(X), B) \rightarrow 
h^G_*(C^\infty_c(X), B)
\end{equation*}
is an isomorphism.
\end{prop}
\proof We use induction on the dimension of $ X $. For $ \dim(X) = 0 $ 
the algebra $ C^\infty_c(X) $ is equipped with the fine bornology and  
$ \theta \Omega_G(C^\infty_c(X) \hat{\otimes} \mathcal{K}_G) $ is 
a projective system of fine spaces. Hence the complexes used in the definition 
of $ HP^G_* $  and $ h^G_* $ are equal and $ v $ is clearly an isomorphism 
in this case. Now assume that the assertion is true for all 
$ G $-simplicial complexes of dimension smaller than $ k $ and that $ \dim(X) = k $. 
Due to proposition \ref{ComplexSplitting} we have a linearly split extension 
of $ G $-algebras of the form 
\begin{equation*} 
 \xymatrix{
     \bigoplus_{j \in J} C^\infty_c(\sigma_j, \partial \sigma_j) \;\; \ar@{>->}[r] &
        C^\infty_c(X) \ar@{->>}[r] \ar@{<.}[l]&  C^\infty_c(X^{k - 1})   
     }
\end{equation*}
where each $ \sigma_j = G/H_j \times \Delta^k $ is an equivariant simplex. 
Using the six-term exact sequences for $ HP^G_* $ 
and $ h^G_* $ obtained from the excision theorem it suffices 
to show that 
\begin{equation*}
v: HP^G_*\Bigl(\bigoplus_{j \in J} C^\infty_c(\sigma_j, \partial \sigma_j), B\Bigr)
\rightarrow 
h^G_*\Bigl(\bigoplus_{j \in J} C^\infty_c(\sigma_j, \partial \sigma_j), B\Bigr)
\end{equation*}
is an isomorphism. Applying excision again we see that in both theories 
$ HP^G_* $ and $ h^G_* $ the $ G $-algebras $ \bigoplus_{j \in J} C^\infty_c(\sigma_j, \partial \sigma_j) $ and 
$ \bigoplus_{j \in J} C_c(G/H_j) $ are equivalent. 
Since $ v $ is compatible with products the assertion follows now 
from the case $ \dim(X) = 0 $ which we have already proved. \qed 
\begin{cor}\label{hcor}
For all finite dimensional $ G $-simplicial complexes $ X $ and $ Y $ we have a
natural isomorphism 
\begin{equation*}
HP^G_*(C^\infty_c(X), C^\infty_c(Y)) \cong
h^G_*(C^\infty_c(X), C^\infty_c(Y)). 
\end{equation*}
This isomorphism is natural with respect to equivariant proper simplicial 
maps in both variables. 
\end{cor}
We come to the next ingredient in the proof of theorem \ref{BSC}. For an arbitrary  $ G $-algebra $ B $ 
one defines a map $ \Tr: \Omega_G(B \cotimes \mathcal{K}_G) \rightarrow \Omega_G(B) $ by
\begin{align*}
\Tr(f(t) &\otimes (x_0 \otimes  k_0) d(x_1 \otimes k_1) \cdots d(x_n \otimes k_n)) \\
&= f(t) \otimes x_0 dx_1 \cdots dx_n \int k_0(r_0,r_1) k_1(r_1,r_2) \cdots k_n(r_n, tr_0) dr_0 \cdots dr_n 
\end{align*}
and
\begin{align*}
\Tr(f(t) &\otimes d(x_1 \otimes k_1) \cdots d(x_n \otimes k_n)) \\
&= f(t) \otimes dx_1 \cdots dx_n \int k_1(r_1,r_2) \cdots k_n(r_n, tr_1) dr_1 \cdots dr_n.
\end{align*}
It is straightforward to check that $ \Tr $ is a map of paramixed complexes. 
\begin{prop}\label{HUnitalProp} Let $ X $ be a locally finite $ G $-simplicial complex. Then the map 
$ \Tr: \Omega_G(C^\infty_c(X) \hat{\otimes} \mathcal{K}_G) \rightarrow \Omega_G(C^\infty_c(X)) $ induces an 
isomorphism on the homology with respect to the equivariant Hochschild boundary.
\end{prop}
\proof Consider the natural commutative diagram 
\begin{equation*} 
    \xymatrix{
    \Omega_G(C^\infty_c(X) \cotimes \mathcal{K}_G) \ar@{->}[r] \ar@{->}[d]^{\Tr} & 
      \Omega_G(C^\infty_c(X)^+ \cotimes \mathcal{K}_G) \ar@{->}[r] \ar@{->}[d]^{\Tr} &
       \Omega_G(\mathbb{C} \cotimes \mathcal{K}_G) \ar@{->}[d] ^{\Tr}\\
 \Omega_G(C^\infty_c(X)) \ar@{->}[r] & 
     \Omega_G(C^\infty_c(X)^+) \ar@{->}[r] &
  \Omega_G(\mathbb{C})  
     }
\end{equation*}
of complexes. According to  proposition \ref{localunits} the algebra $ C^\infty_c(X) $ has local units. 
The same is true for $ C^\infty_c(X) \cotimes \mathcal{K}_G $. Hence the horizontal maps induce long exact sequences in 
homology. Proposition 16.2 in \cite{Voigtperiodic} shows that 
$ \Tr: \Omega_G(B \cotimes \mathcal{K}_G) \rightarrow \Omega_G(B) $ is a linear homotopy equivalence for every 
unital $ G $-algebra $ B $. We deduce that the middle and right vertical arrows induce isomorphisms in homology. 
Consequently, the left vertical arrow induces an isomorphism in homology as well. \qed \\
Let $ X $ be a $ G $-simplicial complex. In section \ref{secHKR} we have studied the equivariant Hochschild-Kostant-Rosenberg 
map $ \Omega_G(C^\infty_c(X)) \rightarrow \mathcal{A}_c(\bar{X}) $. We compose this map with 
the projection onto the elliptic part of $ \mathcal{A}_c(\bar{X}) $ to obtain a map 
$$
\alpha: \Omega_G(C^\infty_c(X)) \rightarrow \mathcal{A}_c(\bar{X})_{\elli} = \mathcal{A}_c(\hat{X}). 
$$
By construction, the hyperbolic part of $ \Omega_G(C^\infty_c(X)) $ is mapped to zero under this map. Let us define a covariant map
$ q: \mathfrak{Fine}(\theta \Omega_G(C^\infty_c(X) \hat{\otimes} \mathcal{K}_G)) \rightarrow \mathcal{A}_c(\hat{X}) $ by composing 
$ \Tr: \Omega_G(C^\infty_c(X) \hat{\otimes} \mathcal{K}_G) \rightarrow \Omega_G(C^\infty_c(X)) $ with the map 
$ \alpha $. \\
Now let $ Q $ be any relatively projective paracomplex of fine covariant pro-modules. Composition with $ q $ yields a map 
\begin{equation*}
f: H_*(\SHom_G(Q, \mathfrak{Fine}(\theta \Omega_G(C^\infty_c(X) \hat{\otimes} \mathcal{K}_G)))) 
\rightarrow H_*(\SHom_G(Q, \mathcal{A}_c(\hat{X}))).
\end{equation*}
Moreover, as explained in section \ref{secBSdef}, we choose a regular projective resolution $ P^\bullet(\hat{X}) $ 
of $ \mathcal{A}_c^\bullet(\hat{X}) $ as above 
and let $ P(\hat{X}) $ be the associated supercomplex. Composition with the chain map 
$ p: P(\hat{X}) \rightarrow \mathcal{A}_c(\hat{X}) $ yields a map
\begin{equation*}
g: H_*(\SHom_G(Q, P(\hat{X}))) \rightarrow H_*(\SHom_G(Q, \mathcal{A}_c(\hat{X})))
\end{equation*}
Remark that the $ \Hom $-complexes occuring in the definition of $ f $ and $ g $ are in fact complexes since all entries 
are paracomplexes of covariant pro-modules.  
\begin{prop}\label{f1f2} Let $ X $ be a finite dimensional locally finite $ G $-simplicial complex. If the action 
of $ G $ on $ X $ is proper the map $ f $ is an isomorphism. 
\end{prop}
\proof We shall treat the elliptic and the hyperbolic parts 
separately. Let us abbreviate $ \Omega_G(X) = \Omega_G(C^\infty_c(X) \hat{\otimes} \mathcal{K}_G) $. \\
First we consider the elliptic part. According to proposition \ref{HUnitalProp} and theorem \ref{HKR} the restriction 
to the elliptic part of the map $ q $ defined above 
induces a quasiisomorphism $ \Omega_G(X)_{\elli} \rightarrow 
\mathcal{A}_c(\hat{X}) $ on the homology with respect to the Hochschild boundary. 
Moreover the natural projection $ E $ on the elliptic part of covariant modules introduced in section \ref{sectotdis}
preserves linearly split exact sequences of covariant modules. In particular 
we obtain a quasiisomorphism $ E\Omega_G(X)_{\elli} \rightarrow E\mathcal{A}_c(\hat{X}) = \mathcal{A}_c(\hat{X}) $ 
which will be denoted by $ q $ again. 
Since $ T = \id $ on $  E\Omega_G(X)_{\elli} $ and $ \mathcal{A}_c(\hat{X}) $ the map 
$ q $ is in fact a map between ordinary mixed complexes. 
Recall from section \ref{secper}that we write $ L^j(M) = F^{j - 1}(M)/F^j(M) $ for the subquotients of the Hodge filtration
of a (para-) mixed complex $ M $. 
The map $ q $ induces chain maps $ L^j E\Omega_G(X)_{\elli} \rightarrow
L^j \mathcal{A}_c(\hat{X}) $ for all $ j $. Since $ q $ is a quasiisomorphism with respect to 
the Hochschild boundary it follows easily that these maps are quasiisomorphism 
of supercomplexes. In particular the corresponding mapping cones are acyclic. \\
We need the following two auxiliary results.
\begin{lemma}\label{BSClemm1} Let $ \phi: D \rightarrow E $ be a morphism of supercomplexes of covariant modules. Assume that 
$ \partial_0: D_0 \rightarrow D_1 $ is zero and 
that $ \partial_1: E_1 \rightarrow E_0 $ is surjective. Then in the mapping 
cone $ C^\phi $ we have $ \im(\partial_1) = E_0 $. Consequently 
the image of $ \partial_1 $ is a direct summand in $ C^\phi_0 = D_1 \oplus E_0 $.
\end{lemma}
\proof By assumption the differential $ \partial_1 $ in $ C^\phi $ has the form 
\begin{equation*}
\begin{pmatrix}
0 & 0 \\
-f & \partial_1 
\end{pmatrix}.
\end{equation*}
Since $ \partial_1: E_1 \rightarrow E_0 $ is surjective the image of 
$ \partial_1 $ in $ C^\phi $ is precisely $ E_0 $. \qed 
\begin{lemma}\label{BSClemm2}
Let $ Q $ be a relatively projective paracomplex of fine covariant pro-modules and let  
$ \phi: Q \rightarrow C $ be a covariant chain map where $ C $ is a constant and acyclic supercomplex. Moreover 
assume that $ C_0 $ admits a direct sum decomposition $ C_0 = K \oplus R $ where 
$ K = \im(\partial_1) = \ker(\partial_0) $. Then $ \phi $ is homotopic to zero. 
Consequently we have $ H_*(\SHom_G(Q,C)) = 0 $. 
\end{lemma}
\proof The map $ \phi_0: Q_0 \rightarrow C_0 = K \oplus R $ may be written as $ \phi_0 = k \oplus r $. 
Since $ \partial_1: C_1 \rightarrow 
K $ is a surjection we find a covariant map $ s: Q_0 \rightarrow C_1 $ such that 
$ \partial_1 s = k $. Hence we may assume without loss of generality that 
$ k = 0 $. Now since $ \phi $ is a chain map and the image of $ \partial_1 \phi_1 $ 
is contained in $ K $ we deduce $ \partial_1 \phi_1 = 0 $. Since 
$ C $ is exact we have $ \im(\phi_1) \subset 
\ker(\partial_1) = \im(\partial_0) $. We may thus construct a 
map $ h: Q_1 \rightarrow C_0 $ such that $ \partial_0 h = \phi_1 $. Furthermore 
we may assume that $ h $ factorizes over $ R $, that is, 
$ h: Q_1 \rightarrow R \rightarrow C_0 $. Hence, up to chain homotopy, the map $ \phi $ satisfies $ \phi_1 = 0 $ and $ k = 0 $ in 
$ \phi_0 = k \oplus r $. Since $ \phi $ is a chain map we now have 
$ 0 = \partial_0 \phi_0 = \partial_0 r $. But $ \partial_0 $ restricted to 
$ R $ is an injection since $ \ker(\partial_0) = K $. This implies $ \phi_0 = 0 $ and 
hence our original map $ \phi $ is homotopic to zero. 
Since we have explicitly shown that any chain map 
$ \phi: Q \rightarrow C $ is homotopic to zero we obtain
$ H_0(\SHom_G(Q,C)) = 0 $. By reindexing $ Q $ we deduce in the same way that
$ H_1(\SHom_G(Q,C)) = 0 $. This finishes the proof. \qed \\
After possibly reindexing, the map 
$$ 
L^j E\Omega_G(X)_{\elli} \rightarrow L^j \mathcal{A}_c(\hat{X}) 
$$ 
satisfies the assumptions of lemma \ref{BSClemm1}. It follows that its mapping cone $ C^j $ satisfies 
the assumptions of lemma \ref{BSClemm2}. 
The short exact sequence 
$$
   \xymatrix{
 L^j \mathcal{A}_c(\hat{X})\;\; \ar@{>->}[r] & 
    C^j \ar@{->>}[r] &  L^j E\Omega_G(X)_{\elli}[1]
     }
$$
of supercomplexes has a covariant splitting and induces  
a short exact sequence 
$$
   \xymatrix{
 \SHom_G(Q_{\elli},L^j \mathcal{A}_c(\hat{X}))\;\; \ar@{>->}[r] & 
    \SHom_G(Q_{\elli},C^j) \ar@{->>}[r] & \SHom_G(Q_{\elli}, L^j E\Omega_G(X)_{\elli}[1])   
     }
$$
of supercomplexes. From lemma \ref{BSClemm2} we deduce 
$ H_*(\SHom_G(Q_{\elli},C^j)) = 0 $. Since the boundary map in this long exact sequence is the map induced by $ q $ 
we obtain an isomorphism 
\begin{equation*}
H_*(\SHom_G(Q_{\elli}, L^j E\Omega_G(X)_{\elli})) \cong H_*(\SHom_G(Q_{\elli}, L^j \mathcal{A}_c(\hat{X})))
\end{equation*}
for all $ j $. \\
Since $ Q $ is assumed to be relatively projective the Hodge filtrations of 
$ E \theta^n \Omega_G(X)_{\elli} $ and $ \theta^n \mathcal{A}_c(\hat{X}) $ induce 
bounded filtrations of the complexes $ \SHom_G(Q_{\elli}, E\theta^n \Omega_G(X)_{\elli}) $ and 
$ \SHom_G(Q_{\elli},\theta^n \mathcal{A}_c(\hat{X})) $, 
respectively. Hence the corresponding spectral sequences converge. The map $ q $ induces a map of spectral sequences which 
gives an isomorphism on the $ E^1 $-terms according to the preceeding discussion. Hence we obtain 
\begin{prop}\label{f1prop} With the notation as above the map $ q $ induces an isomorphism 
\begin{equation*}
H_*(\SHom_G(Q_{\elli}, E\theta^n \Omega_G(X)_{\elli}))
\cong H_*(\SHom_G(Q_{\elli},\theta^n \mathcal{A}_c(\hat{X})))
\end{equation*}
for all $ n $. 
\end{prop}
Let $ M $ be a paramixed complex and set $ M^n = \SHom_G(Q_{\elli}, \theta^n M) $.
Since the structure maps in $ \theta M $ are surjections and 
$ Q_{\elli} $ is relatively projective the structure maps in the inverse system 
$ (M^n)_{n \in\mathbb{N}} $ are surjective. This implies 
$ \varprojlim^1 M^n = 0 $. Therefore we obtain a short exact sequence 
\begin{equation*}
\xymatrix{
    \varprojlim_n  M^n \;\; \ar@{>->}[r] & \prod_{n \in \mathbb{N}} M^n \ar@{->>}[r]^{\id - \sigma} &
        \prod_{n \in \mathbb{N}} M^n }
\end{equation*}
of supercomplexes where $ \sigma $ denotes the structure maps in $ (M_n)_{n \in \mathbb{N}} $. 
This induces a long exact sequence  
\begin{equation*}
\xymatrix{
 {H_0(\SHom_G(Q_{\elli}, \theta M))\;} \ar@{->}[r] \ar@{<-}[d] &
     H_0(\prod_{n \in \mathbb{N}} M^n) \ar@{->}[r] &
       H_0(\prod_{n \in \mathbb{N}} M^n)\ar@{->}[d] \\
   {H_1(\prod_{n \in \mathbb{N}} M^n)\;} \ar@{<-}[r] &
    {H_1(\prod_{n \in \mathbb{N}} M^n)}  \ar@{<-}[r] &
     {H_1(\SHom_G(Q_{\elli}, \theta M))} \\
}
\end{equation*}
in homology. \\
The map $ q $ induces a morphism between these exact sequences for $ M = E\Omega_G(X)_{\elli} $ 
and $ M = \mathcal{A}_c(\hat{X}) $. Hence we obtain an isomorphism  
\begin{equation*}
H_*(\SHom_G(Q_{\elli}, E \theta \Omega_G(X)_{\elli} \cong H_*(\SHom_G(Q_{\elli}, \mathcal{A}_c(\hat{X}))).
\end{equation*}
According to proposition \ref{Econtr} the canonical projection $ \theta \Omega_G(X)_{\elli} \rightarrow E\theta \Omega_G(X)_{\elli} $ 
is a covariant homotopy equivalence. 
Hence we finally conclude that the map
\begin{align*}
f: H_*(\SHom_G(Q_{\elli}, \mathfrak{Fine}(\theta \Omega_G(C^\infty_c(X) \cotimes \mathcal{K}_G))_{\elli})) 
\rightarrow H_*(\SHom_G(Q_{\elli}, \mathcal{A}_c(\hat{X}))_{\elli})) 
\end{align*}
is an isomorphism on the elliptic part. \\
It remains to treat the hyperbolic part. Since $ X $ is proper it follows from the 
equivariant Hochschild-Kostant-Rosenberg theorem 
\ref{HKR} and proposition \ref{HUnitalProp} that the hyperbolic part 
$ \Omega_G(X)_{\hyp} $ of $ \Omega_G(X) $ is acyclic with respect to the Hochschild boundary. 
Hence the associated subquotiens $ L^j\Omega_G(X)_{\hyp} $ of the Hodge filtration 
are covariantly contractible paracomplexes. We can now proceed as above to obtain 
\begin{equation*}
H_*(\SHom_G(Q_{\hyp}, \theta\Omega_G(X)_{\hyp})) = 0. 
\end{equation*}
It follows that the map 
\begin{align*}
f: H_*(\SHom_G&(Q_{\hyp}, \mathfrak{Fine}(\theta\Omega_G(X))_{\hyp})) 
\rightarrow H_*(\SHom_G(Q_{\hyp}, \mathcal{A}_c(\hat{X})_{\hyp})) 
\end{align*}
is an isomorphism because both sides are zero. This shows that $ f $ is an isomorphism on the hyperbolic part. \qed \\
The assumption on $ X $ being proper was only used to prove the isomorphism of the hyperbolic parts in proposition \ref{f1f2}. 
Hence we obtain the following general statement for the elliptic part. 
\begin{prop} \label{f1f2ell} Let $ X $ be any finite dimensional locally finite $ G $-simplicial complex. Then the restriction of the map $ f $ to the 
elliptic part is an isomorphism. 
\end{prop}
Now we study the map $ g $ from above. 
\begin{prop}\label{g1g2} The map $ g $ is an isomorphism for any finite dimensional locally finite $ G $-simplicial complex $ X $. 
\end{prop}
\proof Recall from section \ref{secBSdef} that $ P(\hat{X}) $ can be viewed as an unbounded mixed complex 
with $ b $-boundary equal to zero and $ B $-boundary equal to the 
differential $ \delta $. By definition, the supercomplex $ P(\hat{X}) $ is the projective system of 
supercomplexes $ \xi P(\hat{X}) = (\xi^n P(\hat{X})) $ given by 
\begin{equation*}
\xi^n P(\hat{X}) = P(\hat{X})_{-(n + 1)}/B P(\hat{X})_{-(n + 2)} \oplus \bigoplus_{i = -n}^{n} P(\hat{X})_i \oplus 
B(P(\hat{X})_i).
\end{equation*}
Since $ P(\hat{X}) $ is bounded above in the sense that $ P(\hat{X})_n = 0 $ for 
$ n > D = \dim(\hat{X}) $ we obtain 
\begin{equation*}
\xi^n P(\hat{X}) = P(\hat{X})_{-(n + 1)}/B P(\hat{X})_{-(n + 2)} \oplus \bigoplus_{i = -n}^{D} P(\hat{X})_i
\end{equation*}
for $ n > D $. We define the Hodge filtration $ F^j $ of $ \xi^n P(\hat{X}) $ for $ n > D $ by 
\begin{equation*}
F^j\xi^n P(\hat{X}) = P(\hat{X})_{-(n + 1)}/B(P(\hat{X})_{-(n + 2)}) \oplus \bigoplus_{i = -n}^j P(\hat{X})_i \oplus 
B(P(\hat{X})_j).
\end{equation*}
Hence $ F^j\xi^n P(\hat{X}) $ is a finite increasing filtration such that  
$ F^{-(n + 2)}\xi^n P(\hat{X}) = 0 $ and $ F^n\xi^n P(\hat{X}) = \xi^n P(\hat{X}) $. 
If we proceed in the same way for $ \mathcal{A}_c(\hat{X}) $ we see that the map 
$ p: P^\bullet(\hat{X}) \rightarrow \mathcal{A}^\bullet_c(\hat{X}) $ induces 
chain maps $ \xi^n P(\hat{X}) \rightarrow \xi^n \mathcal{A}_c(\hat{X}) $ which are compatible with the 
filtrations. 
By construction, the map $ p $ is a quasiisomorphism with respect to the boundary 
$ B $. If we denote again by $ L^j $ the subquotients of the Hodge filtration it follows 
that $ p: L^j\xi^n P(\hat{X}) \rightarrow L^j\xi^n \mathcal{A}_c(\hat{X}) $ is a quasiisomorphism 
for each $ j $ and $ n > D $. Hence the corresponding mapping cone $ C^j $ is acyclic. 
Since $ Q $ is relatively projective we see in the same way as in the proof 
of proposition \ref{f1f2} that the map 
\begin{equation*}
H_*(\SHom_G(Q, \xi^n  P(\hat{X})))
\rightarrow H_*(\SHom_G(Q, \xi^n \mathcal{A}_c(\hat{X})))
\end{equation*}
is an isomorphism for $ n > D $. \\
For an unbounded mixed complex $ M $ we set $ M^n = \SHom_G(Q, \xi^n M) $. 
Since for $ M = P(\hat{X}) $ and $ M = \mathcal{A}_c(\hat{X}) $ the projective system 
$ \xi M = (\xi^n M)_{n \in \mathbb{N}} $ 
is isomorphic to the projective systems $ (\xi^n M)_{n > D} $ we obtain 
as in the proof of proposition \ref{f1f2} long exact sequences
\begin{equation*}
\xymatrix{
 {H_0(\SHom_G(Q_{\elli}, \theta M))\;} \ar@{->}[r] \ar@{<-}[d] &
     H_0(\prod_{n \in \mathbb{N}} M^n) \ar@{->}[r] &
       H_0(\prod_{n \in \mathbb{N}} M^n)\ar@{->}[d] \\
   {H_1(\prod_{n \in \mathbb{N}} M^n)\;} \ar@{<-}[r] &
    {H_1(\prod_{n \in \mathbb{N}} M^n)}  \ar@{<-}[r] &
     {H_1(\SHom_G(Q_{\elli}, \theta M))} \\
}.
\end{equation*}
Comparing these exact sequences for $ P(\hat{X}) $ and $ \mathcal{A}_c(\hat{X}) $ 
proves the claim. \qed \\
Propositions \ref{f1f2}, \ref{f1f2ell} and \ref{g1g2} yield the following theorem. 
\begin{theorem}\label{BSC2} Let $ X $ be a finite dimensional locally finite $ G $-simplicial complex. 
If the action of $ G $ on $ X $ is proper the paracomplexes 
$ \mathfrak{Fine}(\theta \Omega_G(C^\infty_c(X) \hat{\otimes} \mathcal{K}_G)) $ and 
$ P(\hat{X}) $ are covariantly homotopy equivalent. \\
The elliptic parts of 
$ \mathfrak{Fine}(\theta \Omega_G(C^\infty_c(X) \hat{\otimes} \mathcal{K}_G)) $ and 
$ P(\hat{X}) $ are covariantly homotopy equivalent even if the action 
of $ G $ on $ X $ is not necessarily proper. 
\end{theorem}
\proof We consider only the case that $ X $ is proper since 
the second assertion is proved in the same way.  
Denote by $ f_1 $ respectively $ f_2 $ the isomorphism $ f $ for 
$ Q = P(\hat{X}) $ and $ Q = \mathfrak{Fine}(\theta \Omega_G(C^\infty_c(X) \hat{\otimes} \mathcal{K}_G)) $. 
Similarly, denote by $ g_1 $ respectivly $ g_2 $ the isomorphism 
$ g $ for $ Q = \mathfrak{Fine}(\theta \Omega_G(C^\infty_c(X) \hat{\otimes} \mathcal{K}_G)) $ 
and $ Q = P(\hat{X}) $. Now let $ x $ be the preimage of $ [p] $ under the isomorphism $ f_1 $ and 
let $ y $ be the preimage of $ [q] $ under the isomorphism $ g_1 $. 
Then we have $ f_1(x) = x \cdot [q] = [p] $ and $ g_1(y) = y \cdot [p] = [q] $. 
Hence $ g_2(x \cdot y) = x \cdot y \cdot[p] = [p] $ and 
$ f_2(y \cdot x) = y \cdot x \cdot [q] = [q] $. 
Since $ g_2 $ and $ f_2 $ are isomorphisms we obtain 
$ x \cdot y = \id $ and $ y \cdot x = \id $. This implies that 
$ \mathfrak{Fine}(\theta \Omega_G(C^\infty_c(X) \hat{\otimes} \mathcal{K}_G)) $ and 
$ P(\hat{X}) $ are covariantly homotopy equivalent. \qed \\
Now we finish the proof of theorem \ref{BSC}. Again we will restrict ourselves to  
the case that $ X $ is a proper $ G $-simplicial complex. \\
Using proposition \ref{BSCProp1} and proposition \ref{f1f2ell} we obtain 
an isomorphism
\begin{equation*}
\bigoplus_{j \in \mathbb{Z}} H^{* + 2j}_G(X,Y) 
\cong H_*(\SHom_G(P(\hat{X}), \mathfrak{Fine}(\theta \Omega_G(C^\infty_c(Y) \hat{\otimes} \mathcal{K}_G)))).
\end{equation*}
Remark that the hyperbolic part of the $ \Hom $-complex vanishes independently
of the fact that the action of $ G $ on $ Y $ may not be proper. 
We apply theorem \ref{BSC2} to deduce 
\begin{align*}
H_*(\SHom_G&(P(\hat{X}), \mathfrak{Fine}(\theta \Omega_G(C^\infty_c(Y) \hat{\otimes} \mathcal{K}_G)))) \\
&\cong H_*(\SHom_G(\mathfrak{Fine}(\theta \Omega_G(C^\infty_c(X) \hat{\otimes} \mathcal{K}_G)), 
\mathfrak{Fine}(\theta \Omega_G(C^\infty_c(Y) \hat{\otimes} \mathcal{K}_G)))).
\end{align*}
Consequently we have an isomorphism
\begin{equation*}
\bigoplus_{j \in \mathbb{Z}} H^{* + 2j}_G(X,Y) 
\cong h_*^G(C^\infty_c(X),C^\infty_c(Y))
\end{equation*}
where $ h_*^G $ denotes the bivariant homology theory introduced before proposition \ref{hprop}.
Combining this with corollary \ref{hcor} we obtain the desired 
identification of equivariant periodic cyclic homology 
with the theory of Baum and Schneider. \\
The naturality of this identification can be checked by straightforward diagram chasing. 
This finishes the proof of theorem \ref{BSC}. \\
Remark that throughout this discussion we used frequently the assumption on $ X $ being a simplicial complex. 
In contrast, it is mainly for convenience to require $ Y $ to be a simplicial complex. 
As soon as an analogue of the equivariant Hochschild-Kostant-Rosenberg theorem for $ Y $ is available,  
the proof presented above can be easily adapted to other classes of spaces. 
We consider briefly the following situation. \\
Let $ G $ be a discrete group and let $ M $ be a $ G $-manifold, that is, a smooth manifold on which $ G $ acts by 
diffeomorphisms. Let $ \mathcal{A}_c(M) $ be the space of smooth differential forms with compact supports
on $ M $ in the usual sense. Then the equivariant Hochschild-Kostant-Rosenberg map $ \alpha: \Omega_G(C^\infty_c(M))_{\elli} 
\rightarrow \mathcal{A}_c(\hat{M}) $ can be defined as in section \ref{secHKR} and one has the following result. 
\begin{theorem}
Let $ G $ be a discrete group and let $ M $ be a $ G $-manifold. 
The equivariant Hochschild-Kostant-Rosenberg map 
$$
\alpha: \Omega_G(C^\infty_c(M))_{\elli} \rightarrow \mathcal{A}_c(\hat{M}) 
$$
induces an isomorphism on homology with respect to the Hochschild 
boundary. 
\end{theorem}
We obtain the following variant of theorem \ref{BSC} in the case of discrete groups. 
\begin{theorem}\label{BSCdis} Let $ G $ be a discrete group and let $ X $ be a proper finite dimensional and locally 
finite $ G $-simplicial complex. Moreover let $ M $ be a $ G $-manifold. Then there exists a 
natural isomorphism
\begin{equation*}
HP^G_*(C^\infty_c(X), C^\infty_c(M)) \cong 
\bigoplus_{j \in \mathbb{Z}} H^{* + 2j}_G(X,M).
\end{equation*}
This isomorphism is natural with respect to equivariant proper simplicial maps in the first variable and 
equivariant proper smooth maps in the second variable. 
\end{theorem}
Let $ Y $ be any locally compact $ G $-space. The equivariant cohomology of $ Y $ with 
$ G $-compact supports is defined by 
$$
H^*(\underline{E}G; Y) = \varinjlim_{K \subset \underline{E}G} H^{* + 2j}_G(K,Y)
$$
where the limit is taken over all $ G $-finite subcomplexes $ K $ of $ \underline{E}G $. Here $ \underline{E}G $ is the universal 
example for proper actions \cite{BCH} which can be chosen to be a simplicial complex. In ~\cite{BS} Baum and Schneider show that 
there is a canonical isomorphism
$$
\bigoplus_{j \in \mathbb{Z}} H^{* + 2j}(\underline{E}G;M) \cong H^*(M,G) 
$$
for all discrete groups $ G $ and $ G $-manifolds $ M $. Here 
$$
H^*(M,G) = \bigoplus_{j \in \mathbb{Z}} H_{2j + *}(EG \times_G T\hat{M}, (EG \times_G T\hat{M}) \setminus \{0\}; \mathbb{C}) 
$$
are the equivariant cohomology groups introduced by Baum and Connes \cite{BC2}. Using theorem \ref{BSCdis} we see that the 
theory of Baum and Connes can be expressed in terms of equivariant cyclic homology. 

\section{Bredon homology and cosheaf homology}\label{secbred}

In this section we review the definitions of equivariant Bredon homology \cite{Bredon1}, \cite{Lueck}
and cosheaf homology \cite{BCH} and compare these theories. Throughout we work with coefficients in the 
complex numbers. \\
Let us begin with Bredon homology. The orbit category $ \Or(G) $ of a totally disconnected group $ G $ has as 
objects all homogenous spaces $ G/H $ where $ H $ is an 
open subgroup of $ G $. The morphisms in $ \Or(G) $ are all $ G $-equivariant continuous maps. One can 
also consider subcategories of $ \Or(G) $ by restricting the class of subgroups. 
We are interested in the class $ \mathcal{F} $ of all compact open subgroups of $ G $. 
The corresponding full subcategory $ \Or(G,\mathcal{F}) $ of $ \Or(G) $ consists of all 
homogeneous spaces $ G/H $ where $ H $ is compact open. \\
If $ \mathcal{C} $ is a small category a covariant (contravariant) $ \mathcal{C} $-vector space is a 
covariant (contravariant) functor from $ \mathcal{C} $ to the category of vector spaces. Morphisms 
of $ \mathcal{C} $-vector spaces are natural transformations. 
More generally one defines covariant and contravariant $ \mathcal{C} $-objects as 
functors with values in arbitrary target categories. If $ G $ is a discrete group 
viewed as a category with one object then a $ G $-vector space is simply a complex 
representation of $ G $. In the sequel we will have to work with certain $ \Or(G,\mathcal{F}) $-vector-spaces 
and $ \Or(G,\mathcal{F}) $-chain complexes. \\
Given a contravariant $ \mathcal{C} $-vector space $ M $ and a covariant $ \mathcal{C} $-vector space 
$ N $ the tensor product $ M \otimes_\mathcal{C} N $ is the direct sum of $ M(c) \otimes N(c) $ over all 
objects $ c \in \mathcal{C} $ divided by all tensor relations $ mf \otimes n - m \otimes fn $ 
for $ m \in M(d), n \in N(c) $ and morphisms $ f: c \rightarrow d $ in $ \mathcal{C} $. \\
Let $ X $ be a proper $ G $-$ CW $-complex. If $ H \subset G $ is an open subgroup 
we have a canonical identification $ X^H = \map_G(G/H,X) $ of the $ H $-fixed point set 
where $ \map_G $ denotes the space of all equivariant continuous 
maps. Using this description of fixed point sets we see that one obtains  
a contravariant functor from $ \Or(G, \mathcal{F}) $ to the category of $ CW $-complexes 
which associates to $ G/H $ the fixed point set $ X^H $. 
Composition with the covariant functor from $ CW $-complexes to chain complexes  
which associates to a $ CW $-complex $ Y $ the cellular chain complex $ C_*(Y) $ with complex coefficients
yields a contravariant $ \Or(G,\mathcal{F}) $-chain complex $ C_*^{\Or(G,\mathcal{F})}(X) $. \\
Next we define a covariant $ \Or(G, \mathcal{F}) $-vector space $ \mathcal{R}_q $ as follows. 
For a compact open subgroup $ H $ of $ G $ set 
$$ 
\mathcal{R}_q(G/H) = K_q(C^*(H)) \otimes_{\mathbb{Z}} \mathbb{C} 
$$
where $ K_* $ denotes topological $ K $-theory and $ C^*(H) $ is the group $ C^* $-algebra of $ H $. Note that 
$ K_0(C^*(H)) = R(H) $ is the representation ring of $ H $ and $ K_1(C^*(H)) = 0 $. 
For a compact totally disconnected group $ H $ the character map induces an isomorphism 
$$ 
K_0(C^*(H)) \otimes_\mathbb{Z} \mathbb{C} \cong \mathcal{R}(H) 
$$ 
where $ \mathcal{R}(H) $ is the ring of conjugation invariant smooth functions on $ H $. \\
We define a chain complex $ C_*^{\Or(G,\mathcal{F})}(X;\mathcal{R}) $ by equipping 
$$
C_*^{\Or(G,\mathcal{F})}(X;\mathcal{R}) = \bigoplus_{p + q = *} C_p^{\Or(G,\mathcal{F})}(X) \otimes_{\Or(G,\mathcal{F})} \mathcal{R}_q
$$
with the differential induced from $ C_*^{\Or(G,\mathcal{F})}(X) $. 
\begin{definition}
Let $ G $ be a totally disconnected group and let $ X $ be a proper $ G $-$ CW $-complex. 
The equivariant Bredon homology of $ X $ (with coefficients in $ \mathcal{R} $) 
is 
$$
\mathcal{B}H^G_*(X) = H_*(C_*^{\Or(G,\mathcal{F})}(X;\mathcal{R})). 
$$ 
\end{definition}
Next we recall the definition of cosheaf homology. Let $ X $ be a simplicial complex. 
We view $ X $ as a category whose objects are the simplices of $ X $ and whose 
morphisms are inclusions of simplices. A cosheaf $ \mathcal{A} $ on $ X $ is a contravariant functor 
from $ X $ to the category of complex vector spaces. More concretely, a cosheaf $ \mathcal{A} $ is specified by vector spaces 
$ \mathcal{A}(\sigma) $ for every simplex $ \sigma \subset X $ and linear maps $ \alpha^\eta_\sigma: \mathcal{A}(\sigma) 
\rightarrow \mathcal{A}(\eta) $ for every inclusion $ \eta \subset \sigma $. These maps are required to 
satisfy $ \alpha_\sigma^\tau = \alpha_\eta^\tau \alpha_\sigma^\eta $ whenever $ \tau \subset \eta \subset \sigma $ 
and $ \alpha_\sigma^\sigma = \id $ for every simplex $ \sigma $. \\
Now let $ G $ be a totally disconnected group and let $ X $ be a proper $ G $-simplicial complex. We are interested 
in the following cosheaf $ \mathcal{R}_X $ on $ X $. For a simplex $ \sigma \subset X $ we define 
$ \mathcal{R}_X(\sigma) = R(G_\sigma) \otimes_{\mathbb{Z}} \mathbb{C} \cong 
\mathcal{R}(G_\sigma) $ where $ G_\sigma $ denotes the stabilizer of $ \sigma $. 
If $ \eta \subset \sigma $ is a face the map $ \iota^\eta_\sigma:\mathcal{R}(G_\sigma) \rightarrow \mathcal{R}(G_\eta) $ 
is given by induction. \\
Let us define a complex $ S_*(X; \mathcal{R}_X) $ as follows. We set 
$$
S_n(X; \mathcal{R}_X) = 
\Biggl(\; \bigoplus_{\dim(\sigma) = n} \mathcal{R}_X(\sigma)\Biggr)/(f[\sigma] - f[-\sigma]\; \text{for}\; f \in \mathcal{R}_X(\sigma))
$$
where the sum is taken over all oriented simplices of $ X $ and $ f[\sigma] $ is our notation for $ f \in \mathcal{R}_X(\sigma) $ viewed as 
element in the summand corresponding to $ \sigma $. Moreover $ -\sigma $ denotes the oriented simplex $ \sigma $ equipped with the opposite 
orientation. 
The boundary $ \partial: S_n(X; \mathcal{R}_X) \rightarrow S_{n - 1}(X; \mathcal{R}_X) $ is defined by 
$$
\partial(f[\sigma]) = \sum_{\begin{smallmatrix} \eta \subset \sigma \\ 
\dim(\eta) = \dim(\sigma) - 1
\end{smallmatrix}}
\iota^\eta_\sigma(f) [\eta]. 
$$
The group $ G $ acts on $ S_*(X; \mathcal{R}_X) $ in a natural way. We let 
$ S^G_*(X; \mathcal{R}_X) $ denote the complex of coinvariants obtained from $ S_*(X; \mathcal{R}_X) $. 
\begin{definition}
Let $ G $ be a totally disconnected group and let $ X $ be a proper $ G $-simplicial complex. 
The cosheaf homology of $ X $ (with values in $ \mathcal{R}_X $) is 
$$
\mathcal{C}H^G_*(X) = H_*(C^G_*(X; \mathcal{R}_X)). 
$$
\end{definition}
Let $ X $ be a $ G $-simplicial complex. Retaining the notation from above we define a map 
$ \phi: S^G_*(X; \mathcal{R}_X) \rightarrow C_*^{\Or(G,\mathcal{F})}(X; \mathcal{R}) $ 
by 
$$
\phi(f[\sigma]) = [\sigma](G_\sigma) \otimes f(G_\sigma) 
$$
where we indicate by the brackets on the right hand side of this equation that the tensor is viewed as an element in 
$ C_*^{\Or(G,\mathcal{F})}(X^{G_\sigma}) \otimes \mathcal{R}(G_\sigma) $ which maps canonically into the 
target of $ \phi $. 
\begin{prop}\label{chambbred}
Let $ G $ be a totally disconnected group and let $ X $ be a proper $ G $-simplicial complex. 
The map $ \phi: S^G_*(X; \mathcal{R}_X) \rightarrow C_*^{\Or(G,\mathcal{F})}(X; \mathcal{R}) $ 
defined above is an isomorphism of chain complexes. 
\end{prop}
\proof First of all one checks that $ \phi $ vanishes on coinvariants and hence yields a well-defined map. 
We compute 
\begin{align*}
\phi \partial(f[\sigma]) = \phi\Biggl(&\sum_{\eta \subset \sigma} \iota^\eta_\sigma(f) \Biggr) 
= \sum_{\eta \subset \sigma} \eta(G^\eta) \otimes \iota^\eta_\sigma(f)(G^\eta) \\
&= \sum_{\eta \subset \sigma} \eta(G^\sigma) \otimes f(G^\sigma)
= \partial([\sigma](G^\sigma) \otimes f(G^\sigma)) = \partial \phi(f[\sigma])
\end{align*}
for a simplex $ \sigma \subset X $ and $ f \in \mathcal{R}(G^\sigma) $. 
Here we use the fact that the simplicial and the cellular chain complex of a simplicial complex 
can be identified. It follows that $ \phi $ is a chain map. 
Let us define a map $ \psi: C_*^{\Or(G,\mathcal{F})}(X; \mathcal{R}) \rightarrow 
S^G_*(X; \mathcal{R}_X) $ by 
$$
\psi([\sigma](H) \otimes f(H)) = \ind_H^{G_\sigma}(f) [\sigma]. 
$$
It is straightforward to check that $ \psi $ is well-defined and by definition we have $ \psi \phi = \id $. 
We calculate 
$$
\phi \psi([\sigma](H) \otimes f(H)) = [\sigma](G^\sigma) \otimes \ind_H^{G^\sigma}(f)(G^\sigma) = [\sigma](H) \otimes f(H) 
$$
and conclude $ \phi \psi = \id $. This finishes the proof. \qed 
\begin{cor}
Let $ G $ be a totally disconnected group and let $ X $ be a proper $ G $-simplicial complex. Then the equivariant Bredon 
homology $ \mathcal{B}H^G_*(X) $ of $ X $ is naturally isomorphic to the cosheaf homology $ \mathcal{C}H^G_*(X) $. 
\end{cor}

\bibliographystyle{plain}

\begin{thebibliography}{99}

\bibitem{BC1} Baum, P., Connes, A., Geometric $ K $-theory for Lie 
groups and foliations, Preprint IHES, 1982
\bibitem{BC2} Baum, P., Connes, A., Chern character for discrete 
groups, in: A f\a^ete of topology, 163 - 232, Academic Press, 1988
\bibitem{BCH} Baum, P., Connes, A., Higson, N., Classifying space for 
proper actions and $ K $-theory of group $ C^* $-algebras, in $ C^* $-algebras: 
1943 - 1993 (San Antonio, TX, 1993), 241 - 291, Contemp. Math. 167, 1994
\bibitem{BS} Baum, P., Schneider, P., Equivariant bivariant Chern character for profinite groups, $ K $-theory 25 (2002), 313 - 353
\bibitem{BL} Bernstein, J., Lunts, V., Equivariant sheaves and functors, 
Lecture Notes in Mathematics 1578, Springer, 1994
\bibitem{BZ} Bernstein, J., Zelevinskii, A., Representations of the group 
$ GL(n,F) $ where $ F $ is a local non-archimedian field, Russian Math. Surveys 31 
(1976), 1 - 68
\bibitem{Bredon1} Bredon, G., Equivariant cohomology theories, Lecture Notes in 
Mathematics 34, Springer, 1967
\bibitem{Bredon2} Bredon, G., Sheaf theory, second edition, Graduate Texts in Mathematics 170, Springer, 1997
\bibitem{Connes1} Connes, A., Noncommutative differential geometry,
Publ. Math. IHES 39 (1985), 257 - 360
\bibitem{Connes2} Connes, A., Noncommutative Geometry, Academic Press,
1994
\bibitem{HN} Higson, N., Nistor, V., Cyclic homology of totally disconnected groups acting on buildings, J. Funct. Anal. 141 
(1996) no. 2, 466 - 495
\bibitem{LueckLNM} L\"uck, W., Transformation groups and algebraic $ K $-theory, Lecture Notes in Mathematics 1408, 
Springer, 1989
\bibitem{Lueck} L\"uck, W., Chern characters for proper equivariant homology 
theories and applications to $ K $- and $ L $-theory, J. Reine Angew. Math. 543 
(2002), 193 - 234
\bibitem{Meyerthesis} Meyer, R., Analytic cyclic cohomology, Preprintreihe SFB 478, Geometrische Strukturen 
in der Mathematik, Münster, 1999
\bibitem{Meyersmoothrep}  Meyer, R., Smooth group representations on bornological vector spaces, Bull. Sci. Math. 128 (2004), 127 - 166 
\bibitem{Spanier} Spanier, E. H., Algebraic Topology, McGraw-Hill, 1966 
\bibitem{Teleman} Teleman, N., Microlocalisation de l'homologie de 
Hochschild, C. R. Acad. Sci. Paris 326 (1998), 1261 - 1264
\bibitem{Voigtthesis} Voigt, C., Equivariant cyclic homology, Preprintreihe SFB 478, Geometrische Strukturen 
in der Mathematik, Münster, 2003
\bibitem{Voigtperiodic} Voigt, C., Equivariant periodic cyclic homology, arXiv:math.KT/0412021 (2004)
\bibitem{Weibel} Weibel, C. A., An introduction to homological algebra, Cambridge Studies in Advanced Mathematics 38, 
Cambridge University Press, 1994
\bibitem{Willis} Willis, G. A., Totally disconnected groups and proofs of 
conjectures of Hofmann and Mukherja, Bull. Austral. Math. Soc. 51 (1995), 489 - 494
\bibitem{Wodzicki1} Wodzicki, M., The long exact sequence in cyclic homology 
associated with an extension of algebras, C. R. Acad. Sci. Paris 306 (1988), 
399 - 403
\bibitem{Wodzicki2} Wodzicki, M., Excision in cyclic homology and
in rational algebraic $ K $-Theory, Ann. of Math. 129 (1989), 591 - 639 

\end{thebibliography}

\end{document}